\documentclass[12pt,twoside,openany,draft]{book}

\usepackage{amsmath}
\usepackage{amsthm}
\usepackage{amssymb}

\def\Lim{\text{Lim}}

\def\raeq{\hbox{$=$\kern-8.8pt\raise4.32pt\hbox{$\scriptstyle\rightarrow$}}}
\def\daeq{\hbox{$=$\kern-8.8pt\raise4.32pt\hbox{$\scriptstyle\rightarrow$}%
 \kern-8.4pt\raise-2.16pt\hbox{$\scriptstyle\rightarrow$}}}

\begin{document}

\begin{center}
Normal Functions and the Problem of the\\
Distinguished Sequences of Ordinal Numbers
\end{center}
\begin{center}
A Translation of
\end{center}
\begin{center}
``Die Normalfunktionen und das Problem der\\
ausgezeichneten Folgen von Ordnungszahlen''\\
by Heinz Bachmann\\
Vierteljahrsschrift der Naturforschenden Gesellschaft in Zurich\\
(www.ngzh.ch)\\
1950(2), 115--147, MR 0036806\\
www.ngzh.ch/archiv/1950\_95/95\_2/95\_14.pdf
\end{center}

Translator's note:
Translated by Martin Dowd, with the assistance of Google translate,
translate.google.com.
Permission to post this translation has been granted by the editors
of the journal.
A typographical error in the original has been corrected,

\begin{center}
\S\ 1 Introduction
\end{center}

1. In this essay we always move within the theory of
Cantor's ordinal numbers.  We use the following notation:

1) Subtraction of ordinal numbers: If $x$ and $y$ are ordinals
with $y\leq x$, let $x-y$ be the ordinal one gets
by subtracting $y$ from the front of $x$, so that
$$y+(x-y)=x$$.

2) Multiplication of ordinals: For any ordinal numbers $x$ and $y$
the product $x\cdot y$ equals $x$ added to itself $y$ times.

3) The numbering of the number classes: The natural
numbers including zero form the first number class, the countably infinite
order types the second number class, etc.; for $k\geq 2$ $\Omega_{k-2}$ is the
initial ordinal of the $k$th number class.  We also use the usual designations
\begin{align*}
  \omega_0&=\omega \\
  \omega_1&=\Omega
\end{align*}

4) The operation of the limit formation: Given a set
of ordinal numbers, the smallest ordinal number $x$,
for which $y \leq x$ for all
ordinal numbers $y$ of this set, is the limit of this set of
ordinals.  If $x_0 \leq x_1 \leq x_2 \dots$ is a monotone,
non-nondecreasing sequence $\{x_\nu\}$
of ordinal numbers $x_\nu$ whose length is a limit ordinal $\lambda$
(i.e., $\nu$ passes through
all ordinal numbers less than $\lambda$, $0\leq \nu<\lambda$), we write for
their limit
$$x=\Lim_{\nu<\lambda}x_\nu$$

5) We make the following three definitions:

 a) If there is an ordinal number $x'$ such that $x=x'+1$,
$x$ is said to be of the first kind.

  If $x\neq 0$ and is not of the first kind, then $x$ is a
limit ordinal.  Then you can let
$x$ be the limit of an ascending sequence of varying length; however,
if $x=\Lim_{\nu<\lambda} x_\nu$ the inequality
$\omega\leq\lambda\leq x$ must hold.

b) $x$ is of the second kind if an ascending sequence of
length $\omega$ can be given whose limit is $x$.

With that, the following alternative definition, which we will often use, is
equivalent:

$x$ is called of the second kind, if an ascending sequence of length
$\lambda$ can be given, whose limit is $x$, and where $\lambda$ is a limit
ordinal of the second number class.

c) $x$ is of the third kind if there is an ascending sequence of
length $\Omega$, whose limit is $x$.
   The last two definitions 5) b) and 5) c) give an exhaustive
and clear classification of the limit ordinals of the third number class.
Because, for such $x$, it cannot hold that
$$x=\Lim_{\mu<\Omega} x_\mu=\Lim_{\nu<\lambda} x_\nu,\>\lambda<\Omega$$
where the two sequences $\{x_\nu\}$ and $\{y_\nu\}$ are ascending.
Were this in fact the case, one could get every ordinal $y_\nu$ as the
smallest $x_{\mu_\nu}$ larger than
$y_\nu$, that is, every ordinal $\nu<\lambda$ is an ordinal $\mu_\nu<\Omega$,
and $\{\mu_\nu\}$ is a monotone,
nondecreasing sequence with $\Lim_{\nu<\lambda}\mu_\nu=\Omega$, but
because of $\lambda<\Omega$ and $\mu_\nu<\Omega$
this can not be.

2. The problem of the distinguished sequences that we here
want to look at, is the following problem: assign to each limit ordinal $y$
of the second
number class a unique ascending sequence $y_0<y_1<y_2<\dots$ of
length $\omega$ such that $y=\Lim_{\mu<\omega}y_\mu$.

This problem is of interest because it is related to the continuum
problem.  Because, if the problem is solved, then you can after the
approach of G.\ H.\ HARDY (1) specify an effective well-order of
an uncountable subset of the continuum by
every ordinal number $y<\Omega$, clearly a number-theoretic problem.
Functions $f_y(x)$ can be assigned by the determinations:
 $f_0(x)=0$.
 If $y=y'+1$ then $f(y(x))=f_{y'}(x)+1$.
 If $y$ is a limit ordinal of the second number class and $\{y_n\}$
  its distinguished sequence let $f_y(x)=f{y_x}(x)$.

The most important approach to solving the problem of distinguished sequences
is that of O. VEBLEN (2).  VEBLEN solves the problem for a large initial
segment of the second number class, by taking a sequence of length
 $\Omega^\Omega+2$ of
normal functions and with its help defined distinguished
sequences for all limit ordinals smaller than a certain one, denoted
$E(1)$, of the second number class.  VEBLEN mentions
the possibility of continuing this procedure (2).

Then CHURCH and KLEENE (3) solve the problem for an even longer initial
segment of the second number class by analyzing some ordinal numbers of this
to assign them a formal representation.  You can, however
show that this method stops somewhere because
a normal function exists whose first critical point
can not be represented.  It must be noticed, that
in principle, such a formal method can not achieve the goal.
You can not assign a formula to every limit ordinal of the second number class,
because every formalism represents only an at most countable set.

3.  In this work we want to
proceed from VEBLEN as discussed in \S 2 and then (\S 4) continue this
(but not to a complete solution to the problem of distinguished
sequences).  Then (\S 5) we give some basic considerations about the
possibility of a complete solution of the problem of the distinguished
sequences, by continuing the procedure of VEBLEN and then (\S 6)
in addition some analogies concerning number-theoretic functions
are considered.  At the end (\S 7) we will relate
the problem of the distinguished
sequences to equivalent other problems related to normal functions.

With VEBLEN we take the following two definitions:

1) A function $\phi(x)$ is called normal if the following
conditions are fulfilled:

a) The set of arguments consists of all ordinal numbers $x$
 with $1<x <\omega_k$,
where $k>0$ is a fixed ordinal number; the value set is a subset
the argument set.

b) Monotonicity: For arbitrary ordinal numbers $x_1$ and $x_2$ of the argument
set with $x_1<x_2$, $\phi(x_1)<\phi(x_2)$.

c) Continuity: If $x$ is a limit ordinal $<\omega_k$, then
$\phi(x)=\Lim_{x'<x}(1+x')$. (We
do not write $\Lim_{x'<x}(x')$ because $\phi(0)$ is undefined.)

For a normal function defined in this way, we to use the term
``normal function of class $k$''; also we define the
monotone increasing number-theoretic functions (they
may be defined as well as the normal functions, but with $k=0$ and under
elimination of the continuity condition) as normal functions of class 0.

For every normal function $\phi$
 $$\phi(x)\geq x$$
for all $x$ of the argument set $(1\leq x<\omega_k)$, further
 $$\phi(x_2)-\phi(x_1)\geq x_2-x_1$$
for any ordinal numbers $x_1,x_2$, with $1\leq x_1<x_2<\omega_k$
 (because $\phi(x_1+\xi)-\phi(x)$
is a normal function of $\xi$; also $\phi(x_1+\xi)-\phi(x_1)\geq\xi$; you put
$\xi=x_2-x_1)$ and the assertion follows), furthermore
  $$\Lim_{x<\omega_k}(1+x)=\omega_k$$

2) When $\phi(x)$ is a normal function of class $k$, the ordinal numbers $x$,
which satisfy the equation $\phi(x)=x$, are the critical points of $\phi$.
VEBLEN proves that there are always critical points, and that the set
of critical points again is the value set of a normal function of class $k$,
the so-called derivative $\phi'$ of $\phi$.

In addition, we define the following: If you have a
normal function $\phi(x)$, we define $\phi^n(x)$ (where $0\leq n<\omega$), the
$n$-fold iteration of $\phi$, as follows:
\begin{align*}
 &\phi^0(x)=x \\
 &\phi^1(x)=\phi(x) \\
 &\phi^{n+1}(x)=\phi(\phi^n(x))
\end{align*}

Let $V\phi$ denote the value set of $\phi$.  The
intersection of sets is indicated by the sign $D$;
e.g.\ $D_{\nu<\lambda}V\phi_\nu$ is the intersection of the value
 sets of a sequence of length
$\lambda$ of normal functions $\phi_\nu$.

\begin{center}
\S\ 2 Three theorems about normal functions
\end{center}

1. In the following normal functions $\phi(x)$ of
class $k$ have value sets consisting of limit ordinals.  We want
in this section to state three theorems
on distinguished sequences that we will use later on.  From the definition
of a normal function follows immediately:

Theorem 1: If $x$ is a limit ordinal $<\omega_k$, and
 $0<x_0<x_1<\dots$ is an ascending
sequence $\{x_\nu\}$ of length $\lambda$, which is a limit ordinal, with
 $x=\Lim_{\nu<\lambda}x_\nu$, then
$\phi(x)=\Lim_{\nu<\lambda}\phi(x_\nu)$,
and the sequence under the limit operator
$\{\phi(x_\nu)\}$ is also an ascending one.

2.  Now it is obvious, to use similar phrases for the values of certain
given normal functions, whose argument is of the first kind, i.e. ascending
sequence for the values $\phi(x+1)$ of a normal function.

We first consider the case that we have a normal function $\phi'$
which is the derivative of a given normal function $\phi$.  VEBLEN (2) shows
 that
one gets all the values of $\phi'$ from the expressions
 $\Lim_{n<\omega}\phi_n(x)$, if for $x$ one uses all arguments;
we can use these expressions for different $x$.
Something more precise is given in VEBLEN;
we only give the following special case:

Theorem 2: The following hold:
\begin{align*}
  &\phi'(1)=\Lim_{n<\omega}\phi^n(1); \\
  &\phi'(x+1)=\Lim_{n<\omega}\phi^n(\phi'(x)+1)
   \hbox{ for all }x\hbox{ with }1<x<\omega_k;
\end{align*}
the sequences under the limit operators are ascending.

Proof: Let's temporarily write $\alpha=\Lim_{n<\omega}\phi^n(\phi'(x)+1)$,
   so for
$1\leq x<\omega_k$
$$\phi(x)<\phi(x)+1\leq\phi(x+1)$$
so (because of the monotonicity of $\phi$, and because the values of $\phi'$
  are critical numbers of $\phi$)
$$\phi(\phi'(x))=\phi'(x)<\phi(\phi'(x)+1)\leq\phi(\phi'(x+1))
   =\phi'(x+1)$$
$$\phi^n(\phi'(x))=\phi'(x)<\phi^n(\phi'(x)+1)\leq\phi^ni(\phi'(x+1))
   =\phi'(x+1)$$
so
$$\phi'(x)<\alpha\leq\phi'(x+1)$$
Since $\alpha$ is also a critical number of $\phi$,
it's in $V\phi'$,
so $\alpha=\phi'(x+1)$.
The sequence for $\alpha$ is in ascending order;
because $V\phi$ only consists of limit numbers,
$$\phi(\phi(x)+1)>\phi(x)+1$$
so
$$\phi^2(\phi'(x)+1)>\phi(\phi'(x)+1)\hbox{ etc}.$$

The proof of the first part of the claim
$\Lim_{n<\omega}\phi^n(1)=\phi'(1)$,
can be done in exactly the same way; you only have
to replace everywhere in the above proof
$0$ for $x$, $\phi(x)$ and $\phi'(x)$.

3.  Now consider, starting from a normal function $\phi$
 of class $k$,
a sequence of some length $\lambda$ (where $\lambda$ is a limit
ordinal $<\omega_k$) of
normal functions $\phi_\mu$ of class $k$ ($\mu$ runs through all
 ordinal numbers less than
$\lambda$) with the property:

I) For every $\mu$ with $1\leq \mu<\lambda$, $V\phi_\mu$ is the intersection
  of all sets of values $V\phi'_\nu$
for the derivatives $\phi'_\nu$ with $0\leq\nu<\mu$.

VEBLEN (2) shows that the enumeration of the common values of these
normal functions is again the value set of a normal function $\psi$ of the same
class:
 $$D_{\mu<\lambda}V\phi_\mu=V\psi$$
In this case the following holds:

Theorem 3: For each ascending sequence $\{\lambda_\nu\}$ whose length
 is a limit ordinal $\pi$,
where $\Lim_{\nu<\pi}\lambda_\nu=\lambda$ (note that $\pi\leq\lambda$),
\begin{align*}
 &\psi(1)=\Lim_{\nu<\pi}\psi_{\lambda_\nu}(1) \\
 &\psi(x+1)=\Lim_{\nu<\pi}\psi_{\lambda_\nu}(\psi(x)+1)
   \hbox{ for all }x\hbox{ with }1<x<\omega_k;
\end{align*}
the sequences under the limit operators are ascending.

\begin{proof}
We temporarily take the abbreviation
$$\alpha_\nu=\phi_{\lambda_\nu}(\psi(x)+1)\hbox{ for }0\leq\nu<\pi$$
Then for $1\leq x<\omega_k$ we have $\psi(x)+1\leq\psi(x+1)$ (because of
monotonicity of $\phi_{\lambda_\nu}$ and because the values
 of $\psi$ are critical points of $\phi_{\lambda_\nu}$).  Also
$$\psi(x)+1\leq\alpha_\nu\leq\phi_{\lambda_\nu}(\psi(x+1))=\psi(x+1)$$
whence
$$\psi(x)+1\leq\Lim_{\nu<\pi}\alpha_\nu\leq\psi(x+1)$$
If $\nu$ is fixed but arbitrary $<\pi$, then all $\alpha_{\nu'}$,
 for $\nu<\nu'<\pi$ lie in $V\phi_{\lambda_\nu}$, so
$\Lim_{\nu<\pi}\alpha_\nu$ is in $V\phi_{\lambda_\nu}$.
 Since this holds for any $\nu<\pi$, $\Lim_{\nu<\pi}\alpha_\nu$ is in $V\psi$;
and so, because of the result obtained above,
 $\Lim_{\nu<\pi}\alpha_\nu=\psi(x+1)$.

The sequence $\{\alpha_\nu\}$ is ascending: Obviously the sequence
 is non-decreasing.
In addition, $\alpha_\nu<\alpha_{\nu+1}$ for $0\leq\nu<\pi$; because of
 condition I) $a_{\nu+1}$ is a
critical point of $\phi_{\lambda_\nu}$
 $$a_{\nu+1}=\phi_{\lambda_\nu}(a_{\nu+1});$$
if $\alpha_\nu=\alpha_{\nu+1}$ then
 $$\phi_{\lambda_\nu}(\psi(x)+1)=\phi_{\lambda_\nu}(\alpha_{\nu+1})$$
would hold, so
 $$\psi(x)+1=\alpha_{\nu+1}$$
would, which is impossible because $\alpha_{\nu+1}$ is a limit ordinal.

The first part of the claim, $\Lim_{\nu<\pi}\phi_{\lambda_\nu}(1)=\psi(1)$,
is proved the same way; one only has to replace the symbols $x$ and $\psi(x)$
 in the above proof by 0.
\end{proof}

4. If we have a sequence of length $\omega_k$ of normal functions $\phi_\mu$
 of class $k$ ($\mu$
goes through all ordinal numbers smaller than $\omega_k$), where for all
 $\mu$ with $1\leq\mu<\omega_k$
the following hold:
\begin{tabular}{rl}
II) a)&$V\phi_\mu$ is a subset of all $V\phi'_\nu$ with $0\leq\nu<\mu$, \\
    b)&If $\mu$ is a limit ordinal $<\omega_k$ then
      $V\phi_\mu=D_{\nu<\mu}V\phi_\nu$,
\end{tabular}\\
then the function formed from the initial numbers
$$\chi(x)=\phi_x(1)$$
is again a normal function of class $k$.  This is proved by VEBLEN (2) for the
special case of the transfinite series of derivatives of a normal function;
but it is valid on the assumption II) in general; because by theorem 3, if
$\mu$ is a limit ordinal $<\omega_k$, $\phi_\mu(1)=\Lim_{\mu<\nu}\psi_\nu(1)$.

On the other hand, for fixed $x>1$ $\phi_\mu(x)$ is not a
 normal function of
 $\mu$; because then
$\phi_\mu(x)>\phi_\mu(1)\geq\mu$, that is, $\phi_\mu(x)>\mu$ for
 $1\leq\mu<\omega_k$.

\begin{center}
\S\ 3 The procedure of VEBLEN in a representation,\\
 which allows generalization
\end{center}

1.  VEBLEN provides a well-ordered, transfinite sequence of length
 $\Omega^\Omega+2$ (that
we denote ${\mathfrak F}_0$) of normal functions $\phi_\eta$ of class 1
($\eta\leq\Omega^\Omega+1$) and defines using these normal functions and using
properties described in the theorems of \S 2 distinguished
sequences for all limit ordinals $y<\phi_{\Omega^\Omega+1}(1)$ ($=E(1)$
 according to the notation
of VEBLEN).  Here, VEBLEN defines for every limit ordinal $<E(1)$
a symbol of the form $\phi(1_1,1_2,\dots,x_\alpha\dots x_\beta)$, divides
 these symbols into seven
classes, and defines for each limit ordinal $<E(1)$, depending on the class
 of the
associated symbol in various ways, a distinguished sequence (2).
These definitions use the theorems of \S 2; their application occurs
but not clearly.  Let's start with the procedure of
VEBLEN, using these symbols in a slightly different representation
to make the application of the theorems of \S 2 more visible
and thus make a generalization of the method possible.  This
generalization will then be carried out in \S 4.

To complete the advertised program, we first need for every
limit ordinal $\eta\leq\Omega^\Omega$ a unique ascending sequence $\{\eta_x\}$
 with length
$\tau_\eta$ (where $\tau_\eta$ is a limit ordinal $\leq\Omega$),
 so that $\eta\raeq\Lim_{x<\tau_\eta}\eta_x$
(the arrow above the equal sign means that $\{\eta_x\}$ is the
unique sequence assigned to $\eta$):

In particular $\eta\raeq\Lim_{\eta'<\eta}(1+\eta')$ for every limit ordinal
$\eta\leq\Omega$ and $\Omega^\Omega\raeq\Lim_{x<\Omega}\Omega^{1+x}$.

If $\eta$ is a limit ordinal of the considered range, then $\eta$ can be
 written as a finite sum
 $$\eta=\sum_{i=0}^n\Omega^{x_i}\cdot y_i+z$$
(4), where $0\leq n<\omega$, $1<y_i<\Omega$,
 $0\leq z<\Omega$ and $\Omega>x_0>x_1>\dots>x_n\geq 1$.
Here we take the following definitions:

If $z$ is of the second kind, let
 $$\eta\raeq\Lim_{z'<z}(\sum_{i=0}^n\Omega^{x_i}\cdot y_i+(1+z'))
   \eqno(\hbox{case }\alpha)$$
Now assume $z=0$, that is
 $$\eta=\sum_{i=0}^n\Omega^{x_i}\cdot y_i$$
If $y_n$ is of the second kind, then let
 $$\eta\raeq\Lim_{y<y_n}(\sum_{i=0}^{n-1}\Omega^{x_i}\cdot y_i+
   \Omega^{x_n}\cdot(1+y))
   \eqno(\hbox{case }\beta)$$
(If $n=0$, the sum $\sum_{i=0}^{n-1}\dots$ must be replaced by 0.)
If $y_n=y_n'+1$, where $y'\geq0$ so that
 $$\eta=\sum_{i=0}^{n-1}\Omega^{x_i}\cdot y_i+
   \Omega^{x_n}\cdot y_n'+
   \Omega^{x_n}$$
you have to distinguish the following sub-cases:

If $x_n$ is of the second kind, let
 $$\eta\raeq\Lim_{x<x_n}(\sum_{i=0}^{n-1}\Omega^{x_i}\cdot y_i+
   \Omega^{x_n}\cdot y_n'+
   \Omega^{1+x})
   \eqno(\hbox{case }\gamma)$$
If $x_n=x_n'+1$, where $x_n\geq 0$, let
 $$\eta\raeq\Lim_{y<\Omega}(\sum_{i=0}^{n-1}\Omega^{x_i}\cdot y_i+
   \Omega^{x_n}\cdot y_n'+
   \Omega^{x_n'}\cdot(1+y))
   \eqno(\hbox{case }\delta)$$
(where $\Omega^0=1$).

Thus, as we can see, the first term $\eta_0$ of each sequence
 $\{\eta_x\}$ is greater than 0; we
note, only sequences with this property will ever be defined
 (cf.\  Note 1 of \S 4).

2. We can now make the statement: It is possible, for every
ordinal number $\eta\leq\Omega^\Omega+1$ to assign a normal function
 $\phi_n$ of class 1,
i.e.\ to construct a sequence ${\mathfrak F}_0$ of length $\Omega^\Omega+2$
  of normal functions of first class,
formed so that the following properties are satisfied for all
  $\eta\leq\Omega^\Omega+1$:

\begin{list}{}{\setlength{\leftmargin}{40pt}\setlength{\labelwidth}{35pt}
 \setlength{\itemsep}{0pt}\setlength{\topsep}{0pt}}
\item[1)]
\begin{list}{}{\setlength{\leftmargin}{15pt}\setlength{\labelwidth}{10pt}
 \setlength{\itemsep}{0pt}\setlength{\topsep}{0pt}}
\item[a)]$\phi_0(x)=\omega^x$.
\item[b)]If $\eta=\eta'+1$ then $\phi_\eta=\phi_{\eta'}'$.
\item[c)]If $\eta$ is of the second kind and $\eta\raeq\Lim_{x<\tau_\eta}\eta_x$
 ($\tau_\eta<\Omega$), then $V\phi_\eta=D_{x<\tau_\eta}V\phi_{\eta_x}$.
 (By specifying $V\phi_\eta$ of course, the function $\phi_\eta$ from the
   enumeration of $V\phi_\eta$, to the
   interval $1\leq x<\omega_k$ is clearly determined.)
\item[d)]If $\eta$ is of the third kind and $\eta\raeq\Lim_{x<\Omega}\eta_x$
   (i.e.\ $\tau_\eta=\Omega$), then $\phi_\eta(x)=\phi_{\eta_x}(1)$.
\end{list}
\item[2)]
The sequences $\{\eta_x\}$ of 1) have the following properties (for all
$x$ with
$0\leq x<\tau_n$, in case 1) d) for all $x$ with $0\leq x<\Omega$):
\begin{list}{}{\setlength{\leftmargin}{15pt}\setlength{\labelwidth}{10pt}
 \setlength{\itemsep}{0pt}\setlength{\topsep}{0pt}}
\item[a)] If $x$ is of the second kind, then $\eta_x$ is of the second kind
  and $\eta_x\raeq\Lim_{x'<x}\eta_{x'}$ (i.e. the
sequence for $\eta_x$ is the a initial subsequence of
that for $\eta$ with $\tau_{\eta_x}=x$).
\item[b)]  If $\eta'$ is a limit ordinal with $\eta_x<\eta'\leq\eta_{x+1}$ then
 for the first term
$\eta_0'$ of the sequence belonging to $\eta'$, $\eta_0'\geq\eta_x$.
\end{list}
\item[3)]
For the sequences $\{\phi_{\eta_x}\}$ the following applies:
\begin{list}{}{\setlength{\leftmargin}{15pt}\setlength{\labelwidth}{0pt}
 \setlength{\itemsep}{0pt}\setlength{\topsep}{0pt}}
\item[]
$\eta_x+1\rightarrow\eta_{x+1}$ for all $x$ with $0\leq x<\tau_\eta$
\item[]
(If $\alpha,\beta$ are ordinal numbers with $\alpha<\beta$ then
 $\alpha\rightarrow\beta$ is an
 abbreviation for: $V\phi_\eta$ is a subset of $V\phi_\alpha$ for all
  $\eta$ with $\alpha<\eta\leq\beta$).
\end{list}
\end{list}

Remark: VEBLEN chooses as initial function $\phi_0(x)=1+x$; we
choose $\phi(x)=\omega^x$, because we want to start from a set of values
 $V\phi_0$
which consists only of limit ordinals, and so from that point on the
functions are normal.  Plus, you want to
choose distinguished sequences
for all ordinals in $V\phi_0$, and for all other limit ordinals of the second
number class.  For $\phi_0(x)=\omega^x$ this can be done,
as we shall see (\S 7).  One could as a starting function
$\phi_0$ take any other class 1 normal function for which this
property is fulfilled.

Proof of the existence of ${\mathfrak F}_0$: The properties 2)
follow from the definitions of the sequences $\{\eta_x\}$,
 as follows.  For $\eta\leq\Omega$ the definition of
the $\phi_\eta$ is made, as in paragraph 3 and 4 of \S 2 (with
 $\omega_k=\Omega$). We
now assume that $\phi_\eta$ is defined for all $\eta<\eta_0$, where
 $\Omega<\eta_0\leq\Omega^\Omega+1$,
so that the conditions 1) 2) 3) are fulfilled for these normal functions, and
show that you can then form $\phi_{\eta_0}$, so that these conditions 
 are satisfied for all $\eta\leq\eta_0$:
\begin{list}{}{\setlength{\leftmargin}{40pt}\setlength{\labelwidth}{35pt}
 \setlength{\itemsep}{0pt}\setlength{\topsep}{0pt}}
\item[a)] If $\eta_0=\eta_0'+1$, then set $\phi_{\eta_0}=\phi_{\eta_0'}'$,
\item[b)] If $\eta_0$ is a limit ordinal and
  $\eta_0\raeq\Lim_{x<\tau_{\eta_0}}\eta_0^{(x)}$, we first prove that
\end{list}
3) is also satisfied for the sequence of functions $\{\phi_{\eta_0}^{(x)}\}$
 belonging to $\eta_0$.  This follows
directly from 2) b), which already holds for the sequence $\{\eta_0(x)\}$:

For a given $x<\tau_{\eta_0}$ we define a finite falling
  sequence $\{\zeta_n\}$:

Let $\zeta_0=\eta_0^{(x+1)}$.

If $\eta_0^{(x)}+1<\zeta_n$ and $\zeta_n=\zeta_n'+1$ let
  $\zeta_{n+1}=\zeta_n'$; then
  $\eta_0^{(x)}+1\leq\zeta_{n+1}<\zeta_n$
and $\zeta_{n+1}\rightarrow\zeta_n$.

If $\eta_0^{(x)}+1<\zeta_n$ and $\zeta_n$ is a limit ordinal with
  $\zeta_n\raeq\Lim_{y<\tau_{\zeta_n}}\zeta_n^{(y)}$
  let $\zeta_{n+1}=\zeta_n^{(1)}$.
By the induction hypothesis
  $\eta_0^{(x)}+1\leq\zeta_{n+1}<\zeta_n$
and $\zeta_{n+1}\rightarrow\zeta_n$ (because $\zeta_n^{(0)}\geq\eta_0^{(x)}$,
  $V_{\phi_\eta}$ is also a subset of $V_{\phi_{\zeta_\eta}(1)}$ for
$\zeta_n^{(1)}<\eta\leq\zeta_n$).

After finitely many regressions there is a natural number $n_0$ with
$\zeta_{n_0}=\eta_0^{(x)}+1$; thus
  $$\eta_0^{(x)}+1=\zeta_{n_0}\rightarrow\zeta_{n_0-1}\rightarrow
    \dots\rightarrow\zeta_1\rightarrow\zeta_0=\eta_0^{(x+1)}$$
so
  $$\eta_0^{(x)}+1\rightarrow\eta_0^{(x+1)}$$
Thus the properties 2) 3) for the sequence $\{\eta_0^{(x)}\}$ and for the
function sequence $\{\phi_{\eta_0}^{(x)}\}$ are proved; it follows that
 this sequence $\{\phi_{\eta_0}^{(x)}\}$ of length
$\tau_{\eta_0}$ of normal functions has the properties I) or II) of \S 2.
One may define:
\begin{align*}
  &V\phi_{\eta_0}=D_{x<\tau_\eta}V\phi_{\eta_0^{(x)}}
   \hbox{ for }\eta_0\hbox{ second kind} \\
  &\phi_{\eta_0}(x)=\phi_{\eta_0}^{(x)}(1)
   \hbox{ for }\eta_0\hbox{ third kind};
\end{align*}
if $\eta_0$ is of the second kind, then $\phi_{\eta_0}(1)\leq\tau_{\eta_0}$.

This proves the existence of ${\mathfrak F}_0$ by means of
 transfinite induction.

We now prove another property of ${\mathfrak F}_0$:

Proposition: For all $\eta<\Omega^\Omega$: If $\eta$ is of the second kind,
 $\eta\raeq\Lim_{x<\tau_\eta}$ and $\phi_\eta(1)=\tau_\eta$
(where $\tau_\eta$ is the length of the sequence belonging to $\eta$)
 there is an
ordinal $\tilde\eta$ of the third kind with $\eta<\eta\tilde\leq\Omega^\Omega$
 and $\tilde\eta\raeq\Lim_{x<\Omega}\tilde\eta_x$, so that
  $$\tilde\eta_x=\eta_x\hbox{ for all }x<\tau_\eta$$
  $$\tilde\eta_{\tau_\eta}=\eta$$ 
$\tau_\eta$ is then a critical point of $\phi_{\tilde\eta}$ i.e.\
 $\tau_\eta$ is in $V\phi_{\tilde\eta+1}$.

Proof: For $\eta<\Omega$ $\tilde\eta=\Omega$.
  If $\eta>\Omega$, then $\eta$ is of second kind, in a case
$\alpha,\beta,\gamma$ (see above). In the case $\alpha$
  $\tilde\eta=\sum_{i=0}^n\Omega^{x_i}\cdot y_i+\Omega$, in case
  $\beta$ $\tilde\eta=\sum_{i=0}^{n-1}\Omega^{x_i}\cdot y_i+\Omega^{x_n+1}$;
in the case of $\gamma$ you have to distinguish two sub-cases: if
  $\eta=\Omega^{x_0}$,
$x_0$ of the second kind, the corresponding $\tilde\eta=\Omega^\Omega$;
 if $\Omega^{x_0}<\eta=\sum_{i=0}^n\Omega^{x_i}\cdot y_i<\Omega^\Omega$, it
follows from property 3) of the $\Omega^\Omega$ distinguished sequence that
 $\phi_{\Omega^{1+x}}$
  $$\phi_\eta(1)>\phi_{\Omega^{x_0}}(1)\geq\tau_{\Omega^{x_0}}
    =x_0\geq x_\eta=\tau_\eta$$
so
  $$\phi_\eta(1)>\tau_\eta$$
If the premises of the proposition are fulfilled, then 1) d)
 $\phi_{\tilde\eta}(\tau_{\tilde\eta})=\phi_{\eta_{\tau_\eta}}(1)=\tau_\eta$,
i.e.\ $\tau_\eta$ is a fixed point of $\phi_{\tilde\eta}$.

3. Using the sequence ${\mathfrak F}_0$ of normal functions
 $\phi_\eta$ of first
class you can get distinguished sequence
$\{y_\eta\}$ of length $\omega$ for
$y<\phi_{\Omega^\Omega+1}(1)=E(1)$,
which we will write with two arrows.

We define $\omega\daeq\Lim_{n<\omega}(1+\eta)$; then assuming that the
 assignment of the
distinguished sequences is carried out for all limit ordinals $y'<y$, where
$\omega<y<E(1)$ is a limit ordinal.  Then you can also for $y$
define the distinguished sequence:

1) If $y$ is not in $V\phi_0$, then we take the unambiguous representation
  $$y=\omega^{x_0}\cdot x_1+x_2$$
where $1\leq x_0<\Omega$, $1\leq x_1<\omega$, $0\leq x_2<\omega^{x_0}$ (2).
 Then $\omega^{x_0}<y$, $x_1<y$,
and $x_2<y$.

If $x_2$ is a limit ordinal and $x_2\daeq\Lim_{n<\omega}x_2^{n)}$,
 then
  $$y\daeq\Lim_{n<\omega}(\omega^{x_0}\cdot x_1+x_2^{(n)})$$

If $x_2=0$ and $x_1=x_1'+1$ ($x_1'\geq 1$), and 
 $\omega^{x_0}\daeq\lim_{n<\omega}\xi_n$, then
  $$y\daeq\Lim_{n<\omega}(\omega^{x_0}\cdot x_1+\xi_n)$$

2) If $y$ is in $V\phi_0$, then there is a last normal function $\phi_\eta$
 with $0\leq\eta\leq\Omega^\Omega$,
so that $y$ is in $V_{\phi_\eta}$, but not in $V_{\phi_{\eta'}}$ for all
 $\eta'$ with $\eta<\eta'\leq\Omega^\Omega+1$. This
can be proved with the same considerations as we (with
more general conditions) give in \S 5; to the
  train of thought
not being interrupted, let us with this reference to \S 5
be content. -- There is thus an $x$ with $1\leq x\leq\Omega$ and
   $$y=\phi_\eta(x)>x$$
If $x$ is a limit ordinal and $x\daeq\Lim(x_n)$, by theorem 1 of \S 2
   $$y\daeq\Lim_{n<\omega}\phi_{\eta}(x_n)$$

If $x=x'+1$ ($x'\geq 0$), one has four cases to distinguish:

a) If $\eta=0$, let
  $$y\daeq\Lim_{n<\omega}(\phi_0(x')\cdot(1+n))$$

b) If $\eta=\eta'+1$, then by the second theorem of \S 2 set
  $$y\daeq\Lim_{n<\omega}\phi_{\eta'}(1)\hbox{ for }x'=0$$
  $$y\daeq\Lim_{n<\omega}\phi_{\eta'}(\phi_\eta(x')+1)\hbox{ for }x'\geq 1$$

c) If $\eta$ is of second kind, then for $x'\geq 1$,
 $y=\phi_\eta(x'+1)>\phi_\eta(x')\geq\phi_\eta(1)\geq\tau_\eta$,
 i.e. $\tau_n<y$;
but also for $x'=0$, $\tau_\eta<y$; because if $y=\phi_\eta(1)=\tau_\eta$,
then, according to the proposition of this paragraph, an $\tilde\eta$ of third
kind would exist.  so that $y$ would be $V\phi_{\tilde\eta}+1$, which is a
contradiction to the fact that $y$ is not in $V\phi_{\eta'}$ for $y'>y$.
If $\eta\daeq\Lim_{x<\tau_\eta}$ and $\tau_\eta\daeq\Lim_{n<\omega}\sigma_n$
 then one can put after theorem 3 of \S 2:
  $$y\daeq\Lim_{n<\omega}\phi_{\eta_{\sigma_n}}(1)\hbox{ for }x'=0$$
  $$y\daeq\Lim_{n<\omega}\phi_{\eta_{\sigma_n}}(\phi_\eta(x')+1)
   \hbox{ for }x'\geq 1$$

d) Finally, if $\eta$ is of the third kind, then
 $y=\phi_\eta(x'+1)=\phi_{\eta'+1}(1)$, and $y$ lies
not in $V\phi_{\eta'}$ for $\eta_{x'+1}<\eta'<\eta$ and $\eta'>\eta$.
 Now $\eta_{x'+1}$ is again third kind. etc.
In any case, there is a finite sequence $\{\zeta_\nu\}$
 with $\eta=\zeta_0>\zeta_1>\dots>\zeta_n\geq 1$,
$1\leq n<\omega$, so that $\zeta_\nu$ is of the third kind for $0\leq\nu<n$
 but $\zeta_n$ not
of the third kind, and $y=\phi_{\zeta_\nu}(1)$ for $0<\nu\leq n$,
 $y=\phi_{\zeta_0}(x'+1)$, $y$ is not in
$V\phi_{\eta'}$ for $\zeta_{\nu+1}<\eta'<\zeta_\nu$ 
 ($0\leq\nu<n$).

If now $\zeta_n=\zeta_n'+1$ then theorem 2 of \S 2 can be applied: Let
  $$y\daeq\Lim_{m<\omega}\phi_{\zeta_n^{\prime m}}(1)$$

If $\zeta_n$ is of the second kind, theorem 3 of \S 2 can be applied.
 This gives then
really a definition for a distinguished sequence for $y$; because it follows
that $y=\phi_{\zeta_n}(1)>\tau_{\zeta_n}$.

First of all, $y\geq\tau_{\eta_n}$ (see the end of the proof of existence of
 ${\mathfrak F}_0$).
If $y=\phi_{\zeta_n}(1)=\tau_{\zeta_n}$ then, according to the further
 proposition given above there would exist an ordinal
$\mu$ of the third kind, with $\zeta_n<\mu\leq\Omega^\Omega$, with $y$
 in $V\phi_{\mu+1}$.
But because $y$ is not in $V\phi_{\eta'}$ for
$\zeta_{\nu+1}<\eta'<\zeta_\nu$ ($0\leq\nu<n$)
 and for $\zeta_n<\eta'\leq\omega^\omega+1$,
  that would be a contradiction.

Now if $\zeta_n\raeq\Lim_{x<\tau_{\zeta_n}}\zeta_n^{(x)}$
 and $\tau_{\zeta_n}\daeq\Lim_{m<\omega}\sigma_m$, then one can set:
  $$y\daeq\Lim_{m<\omega}\phi_{\zeta_n^{(\sigma_m)}}(1)$$

4. We now want to compare this theory with the theory of CANTOR's
 $\epsilon$-numbers
(5.6).  After CANTOR one refers to the critical points of
$\phi_0(x)=\omega^x$ as $\epsilon$-numbers.  One usually designates the first
 $\epsilon$-number $\phi_1(1)$
with $\epsilon$ itself.  These $\epsilon$-numbers play a major role
in the problem of the distinguished sequences:

Just making use of the three operations addition,
 multiplication
and exponentiation of ordinal numbers, and multiple uses
of the constants 1 and $\omega$, all ordinal
 numbers smaller
than $\epsilon$ are obtained; also very easily (in a well-known manner)
 distinguished sequences
may be defined for all limit ordinals of this set.  For $\epsilon$ you need
a new symbol, if one does not want to introduce new operations, and
a new definition, e.g.\
 $\epsilon=\omega^{\omega^{\omega\dots}}\Bigr\}1+n$.
 $\epsilon=\lim_{n<\omega}\omega^{\omega^{\omega\dots}}\Bigr\}1+n$.
After adding this
new symbols, the represented ordinal numbers are those smaller than
the first critical point of the normal function $\epsilon^x$, which we
 call $\epsilon_1$,
and distinguished sequences for such limit ordinals
may be defined; furthermore one defines e.g.\
 $\epsilon_1=\Lim_{n<\omega}\epsilon^{\epsilon^{\epsilon\dots}}\Bigr\}1+n$.
Continuing,
one reaches the first critical point $\epsilon_2$ of $\epsilon_1^x$, etc.
 We still want to arrange
that we only use simple indexes (not indices of indexes).

We therefore obtain a similar procedure to that of VEBLEN;
by diluting the value sets of normal functions
new normal functions are formed: we define for $1\leq\nu<\Omega$
 $\epsilon_{\nu+1}$ as the
first critical point of $\epsilon_\nu^x$; and for
 $\nu=\Lim_{\nu'<\nu}(1+\nu')$,
 $\epsilon_\nu=\Lim_{\nu'<\nu}\epsilon_{1+\nu'}$; with this
we have a normal function $\epsilon_\nu$ (as a function of the index $\nu$).
 In the foregoing
the ordinal numbers using the three operations are needed for the first
critical point of the normal function $\epsilon_\nu$, which we
 denote by $\eta$,
again a new symbol; then the representation succeeds (and also the
 definition
the distinguished sequences) up to the first critical point $\eta_1$ of
 $\eta^x$.
Defining for $1\leq\nu<\Omega$ $\eta_\nu+1$ as the first critical point of
 $\eta_\nu^x$, and for
$\nu=\Lim_{\nu'<\nu}(1+\nu')$,
 $\eta_\nu=\Lim_{\nu'<\nu}\eta_{1+\nu'}$,
we obtain again a normal function $\eta_\nu$
(as a function of $\nu$); to continue you only need a new symbol for
 the first critical point
(e.g.\ $\zeta$), etc.  We want to continue this process until we
have a sequence of length $\omega$ of symbols $1,\omega,\epsilon,\eta\dots$.
 We denote the
limit of this sequence with $\lambda$.  On the basis of the three operations
 addition,
multiplication and exponentiation, using only simple indices and
with finite use of symbols from the above series of length $\omega$.
all ordinal numbers of the second number class, smaller than $\lambda$,
can be represented, and also with the help of the three operations
 distinguished
sequences defined for the limit ordinals less than $\lambda$.

We now want to compare this new sequence of normal functions with 
the sequence ${\mathfrak F}_0$ of VEBLEN:
The derivative of the normal function $\alpha^x$ (where
$\omega<\alpha<\Omega$) is the function $\phi_1(\beta+(x-1))$, where
 $\phi_1(\beta)$ is the first
 $\epsilon$-number
$>\alpha$ (2).  Now, $\epsilon_1$ is the first $\epsilon$-number after
 $\epsilon$, i.e.\ $\epsilon_1=\phi_1(2)$,
 furthermore $\epsilon_{\nu+1}$ is the
first $\epsilon$-number after $\epsilon_\nu$; in general
  $$\epsilon_\nu=\phi(1+\nu)\hbox{ for }1\leq\nu<\Omega.$$
$\eta$ is the smallest ordinal satisfying the relation
 $\eta=\epsilon_\eta=\phi_1(1+\eta)=\phi_1(\eta)$;
thus $\eta=\phi_2(1)$.  Similarly, $\eta_1$ is the first $\epsilon$-number
 after $\eta_1$, so $\eta_1=\phi_1(\eta+1)$,
and $\eta_{\nu+1}$ the first $\epsilon$-number after $\eta_\nu$; in
 general
 $$\eta_\nu=\phi_1(\eta+1)\hbox{ for }1\leq\nu<\Omega.$$

$\zeta$ is the smallest ordinal $>\eta$ satisfying the relation
 $\zeta=\eta_\zeta=\phi_1(\eta+\zeta)=\phi_1(\zeta)$
(using $\eta+\zeta=\zeta$ (5) (6)), so $\zeta=\psi_2(2)$.
 Further $\zeta_\nu=\phi_1(\zeta+\nu)$,
 etc.  Continuing
  $$\lambda=\phi_2(\omega);$$
the procedure of VEBLEN gives thus already in its beginning distinguished
sequences for all limit ordinals $<\lambda$; it immediately delivers much
 stronger thinnings
of the sets of values than the method with the $\epsilon$-numbers
 described above.

\begin{center}
\S\ 4 Continuation of the procedure of VEBLEN
\end{center} 
1. In the last section we first have for every limit ordinal
 $\eta\leq\Omega^\Omega$ a
unique ascending sequence $\{\eta_\nu\}$ of length $\tau_\eta$
with $\eta\raeq\Lim_{x<\tau_\eta}\eta_x$.  To continue the sequence
 ${\mathfrak F}_0$ of
 normal functions $\phi_\eta$ of first class,
one must give to further limit ordinals of the third number class such
ascending sequences.  It is now obvious to perform this task with the help of
  normal functions of
second class, in an analogous way as VEBLEN
with the help of normal functions of first class defines the distinguished
 sequences of
limit ordinals of the second number class.

First, we define a sequence of length $\omega_2+2$ of normal functions
 $F_\zeta(\xi)$ of
second class ($\zeta$ traverses all ordinal numbers $\leq\omega_2+1$) by
the following definitions:

\begin{gather*}
F_0(\xi)=\Omega^\xi \\
F_{\zeta+1}=F_\zeta'\hbox{ for }0\leq\zeta\leq\omega_2 \\
VF_\zeta=D_{\zeta'<\zeta}VF_{\zeta'},\hbox{ if }\zeta%
\hbox{ is a limit ordinal }<\omega_2 \\
F_{\omega_2}(\xi)=F_\xi(1)
\end{gather*}

2.  With the help of this sequence of normal functions $F_\zeta$ of
 second class we want
for every limit ordinal $\eta\leq F_{\omega_2+1}(1)$ a unique ascending
 sequence $\{\eta_x\}$,
whose length is a limit ordinal $\tau_n\leq\Omega$, so that
 $\eta\raeq\Lim_{x<\tau_\eta}\eta_x$.

Here we have to proceed in three steps:

1) First, we give the initial definition:
  $$\eta\raeq\Lim_{\eta'<\eta}(1+\eta'),\hbox{ if }\eta
    \hbox{ is a limit ordinal }\leq\Omega$$

2) In a second step we then assign to each additional limit ordinal
$\eta\leq F_{{\omega_2}+1}(1)$ a function $\Phi_\eta(\xi)$,
a normal function of second class,
so that\\
\hspace*{1in}$\eta=\Phi_\eta(\xi_\eta)$\\
\hspace*{1in}$\xi_\eta$ a certain limit ordinal $<\omega_2$\\
and\\
\hspace*{1in}$\Phi_\eta(\xi)>\xi$ for all $\xi$ with $1\leq\xi\leq\xi_\eta$\\
(note that the ``function'' is a
normal function $\phi_\eta(\xi)$, only defined for $1\leq\xi\leq\xi_\eta$
 and not for $1\leq\xi<\omega_2$ like
this would have to be at a normal second-class function):

a) If $\Omega<\eta\leq F_{\omega_2+1}(1)$ and $\eta$ is not in
  $VF_\omega$, you have for $\eta$ a
clear representation of the form
  $$\eta=\Omega^{\xi_0}\cdot\xi_1+\xi_2$$
where $1\leq\xi_0<\omega_2$, $1\leq\xi_1<\Omega$, $0\leq\xi_2<\Omega^{\xi_0}$.

If $\xi_2$ is a limit ordinal, let
  $$\xi_\eta=\xi_2,\;\Phi_\eta(\xi)=\Omega^{\xi_0}\cdot\xi_1+\xi
   \hbox{ for }1\leq\xi\leq\xi_\eta
   \eqno(Case\ A)$$

If $\xi_2=0$ and $\xi_1=\xi_1'+1$ ($\xi_1'\geq 1$), let
  $$\xi_\eta=\Omega^{\xi_0},\;\Phi_\eta(\xi)=\Omega^{\xi_0}\cdot\xi_1'+\xi
   \hbox{ for }1\leq\xi\leq\xi_\eta
   \eqno(Case\ B)$$

If $\xi_2=0$ and $\xi_1$ is a limit ordinal, let
  $$\xi_\eta=\xi_1,\;\Phi_\eta(\xi)=\Omega^{\xi_0}\cdot\xi
   \hbox{ for }1\leq\xi\leq\xi_\eta
   \eqno(Case\ C)$$

b) If $\eta$ lies in $VF_0$, $\Omega<\eta\leq F_{\omega_2+1}(1)$,
then there is a last
normal function $F_\zeta(\xi)$ with $0\leq\zeta\leq\omega_2+1$,
 such that $\eta$ is in $VF_\zeta$, but not in
 $VF_{\zeta'}$
for all $\zeta'$ with $\zeta<\zeta'\leq\omega_2+1$.  This is by analogy
 to considerations
as we will give in \S 5.  So you have a $\zeta_\eta'$ with
  $1\leq\xi_\eta'<\omega_2$ and
$$\eta=F_\zeta(\xi_\eta')>\xi_\eta'$$

If $\xi_\eta'=\xi_\eta''+1$ for $\xi_\eta''\geq1$, then $\zeta<\omega_2$, 
 otherwise $\eta$ would be in $VF_{\omega_2}$; one has
the following cases:

$\zeta=0,\,\xi''_\eta\geq 1$: Then\\
\hspace*{.75in}$\xi_\eta=\Omega,\,\Phi_\eta(\xi)=\Omega^{\xi_\eta''}\cdot\xi
 \hbox{ for }1\leq\xi\leq\xi_\eta$\hfill($Case\ D$)

$\zeta=\zeta'+1,\,\zeta'<\omega_2,\,\xi_\eta''\geq 1$: Then\\
\hspace*{.75in}$\xi_\eta=\omega,\,\Phi_\eta(\xi)=
 F_{\zeta'}^\xi(F_\zeta(\xi_\eta'')+1)
 \hbox{ for }1\leq\xi\leq\xi_\eta,\,\Phi_\eta(\xi_\eta)=\eta$\hfill($Case\ E$)

\noindent The second theorem of \S 2 applies.

$\zeta$ a limit ordinal of the third number class,
 $\xi_\eta''\geq 1$: Then let\\
\hspace*{.75in}$\xi_\eta=\zeta,\,\Phi_\eta(\xi)=
 F_\xi(F_\eta(\xi_\eta'')+1)
 \hbox{ for }1\leq\xi<\xi_\eta,\,\Phi_\eta(\xi_\eta)=\eta$\hfill($Case\ F$)

\noindent Here $\Phi_\eta(\xi)$ really is a normal function; because if
 $\zeta'$ is a limit ordinal $<\zeta$
then (by theorem 3 of \S 2, and because $F_\zeta(\xi_\eta'')$ is a critical
 point of $F_{\zeta'}$)
  $$\Phi_\eta(\zeta')=F_{\zeta'}(F_\zeta(\xi_\eta'')+1)=
   \Lim_{\zeta''<\zeta'}F_{1+\zeta''}(F_{\zeta'}(F_\zeta(\xi_\eta''))+1)$$
  $$=\Lim_{\zeta''<\zeta'}F_{1+\zeta''}F_\zeta(\xi_\eta'')+1)=
   \Lim_{\zeta''<\zeta'}\Phi_\eta(1+\zeta'')$$
furthermore
  $$\Phi_\eta(\xi)>F_\xi(1)\geq\xi\hbox{ for }1\leq\xi\leq\xi_\eta$$
In case F, theorem 3 of \S 2 is applied.

$\zeta=\omega_2,\,\xi_\eta''\geq 0$:
 Then $\eta=F_{\xi_\eta''+1}(1)$ and\\
\hspace*{.75in}$\xi_\eta=\omega,\,\Phi_\eta(\xi)=
 F_{\xi_\eta''}^\xi(1)
 \hbox{ for }1\leq\xi<\xi_\eta,\,\Phi_\eta(\xi_\eta)=\eta$\hfill($Case\ E'$)\\
The second theorem of \S 2 is also used here.

$\zeta=\omega_2+1,\,\xi_\eta''=0,\,\eta=F_{\omega_2+1}(1)$: Then\\
\hspace*{.75in}$\xi_\eta=\omega,\,\Phi_\eta(\xi)=
 F_{\omega_2}^\xi(1)
 \hbox{ for }1\leq\xi<\xi_\eta,\,\Phi_\eta(\xi_\eta)=\eta$\\
This case is also taken under Case E${}'$.

If $\xi_\eta'$ is a limit ordinal, $\eta=F_\zeta(\xi_\eta')>\xi_\eta'$,
 then if $\eta$ is greater than the first critical point of $F_\zeta$,
then there exists is a final critical point $\rho$ of $F_\zeta$, such that
$\rho=F_\zeta(\rho)=F_{\zeta+1}(\mu)$, $\mu\geq 1$,
 $F_{\zeta+1}(\mu)<\eta<F_{\zeta+1}(\mu+1)$.
This follows because,
 $\xi_\eta'$ is not a critical point of $F_\zeta$, and the critical
 points of $F_\zeta$ are the
value set of a normal function.  If $\eta$ is smaller than the first
 critical point of $F_\zeta$, we set $\rho=0$.  Now set\\
\hspace*{.75in}$\xi_\eta=\xi_\eta'-\rho\,\Phi_\eta(\xi)=F_\zeta(\rho+\xi)
  \hbox{ for }1\leq\xi\leq\xi_\eta$\hfill($Case\ G$)

The normal functions $\Phi_\eta(\xi)$ thus defined satisfy in all cases A to G
the required conditions.

3) In a third step we can now define the sequences $\{\eta_x\}$ for all
 limit ordinals $\eta\leq F_{{\omega_2}+1}(1)$:
 in the first step this was done for $\eta\leq\Omega$.
Suppose the definitions are given for all limit ordinals
 $\eta'<\eta$, where $\eta$
is a limit ordinal with $\Omega<\eta\leq F_{{\omega_2}+1}(1)$.
Then you can also
define such a sequence for $\eta$; because after the second step you have a
clearly defined normal function $\phi_\eta(\xi)$ with the above described
properties, so that
  $$\eta=\Phi_\eta(\xi_\eta),\,\xi_\eta\hbox{ a limit ordinal }{<}\eta$$
If $\xi_\eta\raeq\Lim_{x<\tau_{\xi_\eta}}\xi_\eta^{(x)}$, then set
 $\eta\raeq\Lim_{x<\tau_{\xi_\eta}}\Phi(\xi_\eta^{(x)})$.
The task of determining
the sequence $\{\eta_x\}$ is thus accomplished
by means of the principle of transfinite induction.
Note that $\tau_\eta=\tau_{\xi_\eta}$
for all $\eta$ with $\Omega<\eta\leq F_{{\omega_2}+1}(1)$.

3. According to the above explanations, it is now possible to create a
 sequence ${\mathfrak F}_1$ of
length $F_{{\omega_2}+1}(1)+1$ of normal functions $\phi_\eta(x)$ of
 first class,
such that for all $\eta\leq F_{{\omega_2}+1}(1)$:

\begin{list}{}{\setlength{\leftmargin}{40pt}\setlength{\labelwidth}{35pt}
 \setlength{\itemsep}{0pt}\setlength{\topsep}{0pt}}
\item[1)]
\begin{list}{}{\setlength{\leftmargin}{15pt}\setlength{\labelwidth}{10pt}
 \setlength{\itemsep}{0pt}\setlength{\topsep}{0pt}}
\item[a)]$\phi_0(x)=\omega^x$.
\item[b)]If $\eta=\eta'+1$ then $\phi_\eta=\phi_{\eta'}'$.
\item[c)]If $\eta$ is of the second kind and $\eta\raeq\Lim_{x<\tau_\eta}\eta_x$
 then $V\phi_\eta=D_{x<\tau_\eta}V\phi_{\eta_x}$.
\item[d)]If $\eta$ is of the third kind and $\eta\raeq\Lim_{x<\Omega}\eta_x$
  then $\phi_\eta(x)=\phi_{\eta_x}(1)$.
\end{list}
\item[2)]
\begin{list}{}{\setlength{\leftmargin}{15pt}\setlength{\labelwidth}{10pt}
 \setlength{\itemsep}{0pt}\setlength{\topsep}{0pt}}
\item[a)] If $\eta\raeq\Lim_{x<\tau_\eta}\eta_x$
and $\mu$ is a limit ordinal $<\tau_\eta$, then
$\eta_\mu\raeq\Lim_{x<\mu}\eta_x$.
\item[b)] If $\eta\raeq\Lim_{x<\tau_\eta}\eta_x$ then
$\eta_x+1\rightarrow\eta_{x+1}$ for all $x$ with $0\leq x<\tau_\eta$.\\
(Again we use the abbreviation $\alpha\rightarrow\beta$
 for: $V\phi_\eta$ is a subset of $V\phi_\alpha$
 all $\eta$ with $\alpha<\eta\leq\beta$.)
\end{list}
\end{list}
We now proceed to prove this claim.  For the beginning of the
sequence ${\mathfrak F}_1$ the above properties are obviously fulfilled,
 because ${\mathfrak F}_0$ is an
initial segment of ${\mathfrak F}_1$ because the sequences
$\{\eta_x\}$ for $\eta\leq\Omega^\Omega$
for the former and in this
paragraph are the same.

To prove the existence of the whole sequence ${\mathfrak F}_1$
we assume that we have
already set up all $\phi_\eta$ for all $\eta<\eta_0$, where
 $\Omega<\eta_0\leq F_{{\omega_2}+1}(1)$, so that the
above conditions 1) 2) are satisfied.  Then we show that $\phi_{\eta_0}$
can be defined, so that the conditions apply to all $\eta\leq\eta_0$.

If $\eta_0$ is a limit ordinal, and condition 2) holds for the sequence
 $\{\eta_0^{(x)}\}$
for $\eta_0$ (so that
 $\eta_0\raeq\Lim_{x<\tau_{\eta_0}}\eta_0^{(x)}$)
then one can, since then the sequence $\{\phi_{\eta_0^{(x)}}\}$ of
normal functions has the property I) or II) of \S 2, for $\eta_0$ of
 the second kind set
$V\phi_{\eta_0}=D_{x<\tau_{\eta_0}}V\phi_{\eta_0^{(x)}}$,
and for $\eta_0$ of the third kind set
 $\phi_{\eta_0}(x)=\phi_{\eta_0^{(x)}}(1)$.
So we just have to prove condition 2) for the sequence $\{\eta_0^{(x)}\}$.

It is very easy to show that 2) a) holds: If $\eta$ is a limit ordinal with
$\Omega<\eta\leq F_{{\omega_2}+1}(1)$, then, according to no.\ 2 of this
 section, one has a limit ordinal $\xi_\eta$, so
that $\xi_\eta=\Phi_{\xi_\eta}(\xi_{\xi_\eta})>\xi_{\xi_\eta}$,
 thus (if $\xi_\eta>\Omega$) a limit ordinal $\xi_{\xi_\eta}$, such that
 $\xi_\eta=\Phi_{\xi_\eta}>\xi_{\xi_\eta}$
etc.  After finitely many steps you get the form
 $$\eta=\Phi_{\zeta_0}\Phi_{\zeta_1}\dots\Phi_{\zeta_n}(\mu_\eta))\dots)$$
(with superfluous parentheses omitted) with
  $$\eta=\zeta_0>\zeta_1>\dots>\zeta_n>\mu_\eta,\,0\leq n<\omega$$
$\Psi_\eta(x)$ (defined at least for $1\leq x\leq\mu_\eta$) is an initial
 segment of a normal function;
because a normal function of a normal function is once again
a normal function.  From the definitions of the sequences $\{\eta_x\}$
 follows, how to
easily consider,
$$\eta_x=\Psi_\eta(1+x);$$
if $\mu$ is a limit ordinal $\leq\mu_\eta$
then (letting $\mu=\Lim_{\mu'<\mu}(1+\mu')$,
$\Phi_{\zeta_n}(\mu)\raeq\Lim_{\mu'<\mu}\Phi_{\zeta_n}(1+\mu')$,
$\Phi_{\zeta_{n-1}}\Phi_{\zeta_n}(\mu)\raeq
 \Lim_{\mu'<\mu}\Phi_{\zeta_{n-1}}\Phi_{\zeta_n}\Phi_{\zeta_n}(1+\mu')$,
 etc.)
$$\Psi_\eta(\mu)\raeq\Lim_{x<\mu}\Psi_\eta(1+x)$$
Further, $\mu_\eta=\tau_\eta$, where $\tau_\eta$ is the length of
 the sequence belonging to $\eta$. -- 
This proves 2) a).

To prove 2) b), we need the following three propositions:

Proposition 1:
Assume that an initial subsequence of ${\mathfrak F}_1$, satisfying necessary
requirements (for $\eta<\eta_0$) is given.
Suppose $F(\xi)$ is a
normal function of second
class with the properties:

$F(\xi)+1\rightarrow F(\xi+1)$ for all $\xi\geq 1$
 with $F(\xi+1)<\eta_0$

$F(\eta)\raeq\Lim_{x<\tau_\eta} F(\eta_x)$,
 if $\eta\raeq\Lim_{x<\tau_\eta} \eta_x$
 for all limit ordinals $\eta$ with $F(\eta)<\eta_0$\\
Let $\alpha$, $\beta$ be ordinal numbers with $1\leq\alpha<\beta$,
 $\alpha\rightarrow\beta$, and $F(\beta)<\eta_0$. -- Then
 $F(\alpha)\rightarrow F(\beta)$.

Proof: a) First, we show by transfinite induction that
$V\phi_{F(\zeta)}$ is a subset of $V\phi_{F(\alpha)}$
for $\alpha<\zeta\leq\beta$: This is because of
 $F(\alpha)+1\rightarrow F(\alpha+1)$
for $\zeta=\alpha+1$.  Now let $\alpha+1<\zeta\leq\beta$
 and let $V\phi_{F(\zeta')}$ be a subset
of $V\phi_{F(\zeta')}$ for all $\zeta'$ with $\alpha<\zeta'<\zeta$;
 we show that then
$V\phi_{F(\zeta)}$ is a subset
of $V\phi_{F(\alpha)}$:

If $\zeta=\zeta'+1$ ($\zeta'>\alpha$), then the claim is proven  because
 of $F(\zeta')+1\rightarrow F(\zeta'+1)$.

If $\zeta$ is second kind and $\zeta\raeq\Lim_{x<\tau_\zeta}\zeta_x$,
 then $F(\zeta)\raeq\Lim_{x<\tau_\zeta}F(\zeta_x)$;
 there is a first $x_0$
such that $F(\zeta_{x_0})>F(\alpha)$; so, by the induction hypothesis
 $V\phi_{F(\zeta_{x_0})}$
is a subset of $V\phi_{F(\alpha)}$, and because of property
 2) b) of the sequence $V\phi_{F(\zeta_x)}$ and
by the definition of $\phi_{F(\zeta)}$, $V\phi_{F(\zeta)}$ is a subset of
 $V\phi_{F(\zeta_{x_0})}$, that is
$V\phi_{F(\zeta)}$ is a subset of
 $V\phi_{F(\alpha)}$.

If $\zeta$ is third kind and $\zeta\raeq\Lim_{x<\Omega}\zeta_x$,
because of property 2) b) of the sequence
$\{\phi_{F(\zeta_x)}\}$
and according to the definition of $\phi_{F(\zeta)}$ $V\phi_{F(\zeta)}$
 is a subset of $V\phi_{F(\zeta_1)}$.
Now, in addition, $\alpha\leq\zeta_1$; for if $\alpha>\zeta_1$,
 then because of property 2) b) of the sequence $\{\phi_{\zeta_x}\}$
$\zeta_1+1\rightarrow\alpha$ would hold, 
 so $\phi_\zeta(1)=\phi_{\zeta_1}(1)<\phi_\alpha(1)$.
 Because of the requirement $\alpha\rightarrow\beta$,
 $\alpha\rightarrow\zeta$,
that is, $\phi_\alpha(1)\leq\phi_\zeta(1)$, which contradicts
previous results.  So $\alpha\leq\zeta_1$ and thus $V\phi_{F(\xi)}$ is a
 subset of $V\phi_{F(\alpha)}$.

b) If $F(\alpha)<\zeta<F(\beta)$, but $\zeta$ is not in $VF$ and
 $F(\gamma)<\zeta<F(\gamma+1)$,
$\alpha\leq\gamma<\beta$, then because of the premise
$F(\gamma)\rightarrow F(\gamma+1)$ $V\phi_\zeta$ is a
subset of $V\phi_{F(\gamma)}$, hence of $V\phi_{F(\alpha)}$.

The first proposition is also ``concrete'';
because one can understand the formation of all normal functions
 $\phi_{\xi'}$,
for $\xi'<\xi$ as $\xi$-fold iteration of the formation of the
 derivative,
and
the formation of all $\phi_{\xi'}$ for $\xi'<F(\xi)$ as $\xi$-fold
 iteration of the operation,
each of $\phi_{F(\xi')}$ over all
  $\phi_\eta$ ($F(\xi')<\eta<F(\xi+1)$)
  on $\phi_{F(\xi'+1)}$.

These two transfinite iterations are to some degree
 similar, since if $\eta$ is a limit,
the $\eta$-fold iterations in both cases exactly match
(because $F(\eta)\raeq\Lim_{x<\tau_\eta}F(\eta_x)$
  when $\eta\raeq\Lim_{x<\tau_\eta}\eta_x$).

Proposition 2: Every ordinal number of $VF_0$ (and also of $V\phi_0$) has the
property that it is equal to all its residues, i.e.\ if $\alpha$
 is in $VF_0$ or $V\phi_0$
and if $\beta<\alpha$, then $\alpha-\beta=\alpha$, or $\beta+\alpha=\alpha$.

Proof: We prove with a slight generalization a fact from
HAUSDORFF (4), that for $\xi\geq 1$ each power $\rho^\xi$ of
the ordinal number $\rho$, which is equal to all its residues, also
 has this property:

Let $\beta<\rho^\xi$, so $\beta=\rho^\eta\cdot\zeta+\mu$, where
$$0\leq\eta<\xi,\,0\leq\zeta<\rho,\,0\leq\mu<\rho^\eta$$
If one defines $\sigma$ by $\rho^\xi=\rho^{\eta+1}+\sigma$, then because
 $\rho\neq 1$
$$\beta+\rho^\xi=\rho^\eta\cdot\zeta_\mu+\rho^\xi\leq\rho^\eta\cdot(\zeta+1)
  +\rho^{eta+1}+\sigma=
  \rho^\eta\cdot(\zeta+1+\rho)+\sigma=\rho^{\eta+1}+\sigma=\rho^\xi$$
If we take for $\rho$ the ordinals $\omega$ and $\Omega$, then we get
 the statement of
second proposition.

Proposition 3: Assuming that all $\phi_\eta$ for $\eta<\eta_0$ satisfy the
requirements 1) 2) (No.\ 3 of \S 4), then
$$\Phi_\eta(\xi)+1\rightarrow\Phi_\eta(\xi+1)$$
for any function $\Phi_\eta$ defined in no.\ 2 of this section, for
$$1\leq\xi<\xi_\eta\hbox{ and }\Phi_\eta(\xi+1)<\eta$$

Proof: For $\Phi_\eta(\xi)$ cases A to G are considered (cf.\ Nr.\ 2
of \S 4).

For cases A and B, the claim is trivial.

In cases C and D, $\phi_\eta(\xi)=\Omega^{\xi_0}\cdot\xi$; because of
 $1\rightarrow\Omega^{\xi_0}$ by the
first proposition (where the function $F$,
 $F(\zeta)=\Omega^{\xi_0}\cdot\xi+\zeta$
is defined)
$$\Phi_\eta(\xi)+1= \Omega^{\xi_0}\cdot\xi+1\rightarrow
  \Omega^{\xi_0}\cdot\xi+\Omega^{\xi_0}= \Phi_\eta(\xi+1)$$
if $\Phi_\eta(\xi+1)<\eta_0$.

Case G: We have $\Phi_\eta(\xi)=F_\zeta(\rho+\xi)$, $\rho\geq 0$.
 We prove that
$F_\zeta(\xi)+1\rightarrow F_\zeta(\xi+1)$ holds for all $\zeta$ with
 $0\leq\zeta\leq\omega_2+1$, provided $F_\zeta(\xi+1)<\eta_0$:

This assertion holds for $\zeta=0$, because as in cases C and D,
 $$F_0(\xi)+1=\Omega^\xi+1\rightarrow\Omega^\xi\cdot 2$$
further
 $$\Omega^\xi\cdot 2 \rightarrow \Omega^\xi\cdot \Omega=F_0(\xi+1)$$
from property 2) b) of the sequence $\{\phi_{\Omega^\xi\cdot x}\}$ and
 from the definition of $\phi_{\Omega^\xi\cdot\Omega}$.
So you have $F_0(\xi)+1\rightarrow F_0(\xi+1)$.

Now suppose the assertion holds for all $F_{\xi'}$ with $\zeta'<\zeta$, where
$0<\zeta\leq\omega_2+1$; we will show that
 it holds for $F_\zeta$:

If $\zeta=\zeta'+1$ then
 $F_\zeta(\xi+1)\raeq\Lim_{n<\omega}F_\zeta^{1+n}(F_\zeta(\xi)+1)$.
 Because of property
2) b) the sequence $\{\phi_{F_{\zeta'}^{1+n}(F_\zeta(\xi)+1)}\}$
 satisfies, using the induction hypothesis,
$$F_{\zeta'}(F_\zeta(\xi))+1= F_\zeta(\xi)+1\rightarrow
  F_{\zeta'}(F_\zeta(\xi)+1)\rightarrow F_\zeta(\xi+1)$$

If $\zeta$ is a limit ordinal $<\omega_2$ with
 $\zeta\raeq\Lim_{x<\tau_\zeta}\zeta_x$ then
  $$F_\zeta(\xi+1)\raeq
    \Lim_{x<\tau_\zeta}F_{\zeta_x}(F_\zeta(\xi)+1)$$
By arguments similar to ones given above
$$F_{\zeta_0}(F_\zeta(\xi))+1=
  F_\zeta(\xi)+1\rightarrow
  F_{\zeta_0}(F_\zeta(\xi)+1)\rightarrow
  F_\zeta(\xi+1)$$

If $\zeta=\omega_2$, then,
 $F_\zeta(\xi+1)= F_{\xi+1}(1)\raeq \Lim_{\eta<\omega}F_\xi^{1+n}(1)$,
so because of property 2) b) for the sequence
 $\{\phi_{x_\xi^{1+n}(1)}\}$
$$F_\zeta(\xi)+1= F_\xi(1)+1\rightarrow F_\zeta(\xi+1)$$

Case E: Then $\Phi_\eta(n)=F_{\zeta^{\prime n}}(\delta)$, where $n<\omega$,
$\xi_\eta''\geq 1$ and
  $$\delta=F_\zeta(\xi_\eta'')+1$$
If $F_{\zeta^{\prime 2}}(\delta)<\eta_0$, one has because of the result
 obtained in case G
 $$F_{\zeta'}(\delta)+1\rightarrow F_{\zeta'}(\delta+1);$$
you use Proposition 1 (with
 $F(\xi)=F_{\zeta'}(F_\zeta(\xi_\eta'')+\xi)$\footnote{The application of
the first proposition is permitted because the first critical
point $F_{\zeta'}(\xi)$ after $F_\zeta(\xi_\eta'')$
is $F_\zeta(\xi_\eta''+1)$, that is the condition
$F(\xi)\raeq\Lim_{x<\tau_\xi}F(\xi_x)$
is always satisfied for the $\xi$ in question.}) 
and Proposition 2,
so because $2\rightarrow F_{\zeta'}(\delta)$
(cf.\ the remark at the end of this proof)
and $F_\zeta(\xi_\eta'')<R_{\zeta'}(\delta)$
$$F_{\zeta'}(\delta+1)=F_{\zeta'}(F_\zeta(\xi'')+2)=F(2)\rightarrow
  F(F_{\zeta'}(\delta))=$$
$$F_{\zeta'}(F_\zeta(\xi_\eta'')+F_{\zeta'}(\delta))=
  F_{\zeta'}(F_{\zeta'}(\delta))=F_{\zeta^{\prime 2}}(\delta)$$
therefore
$$F_{\zeta'}(\delta)+1\rightarrow F_{\zeta^{\prime 2}}(\delta)$$

In an analogous manner, for $1\leq n<\omega$ and
 $F_{\zeta'}^{n+2}(\delta)<\eta_0$ from
$F_{\zeta'}^n(\delta)+1\rightarrow F_{\zeta'}^{n+1}(\delta)$
because of Proposition 1 (with
 $F(\xi)=F_{\zeta'}(F_\zeta(\xi_\eta'')+\xi)$ and
Proposition 2:

$$F_{\zeta'}^{n+1}(\delta)+1=F_{\zeta'}F_{\zeta'}^n(\delta))+1\rightarrow
  F_{\zeta'}(F_{\zeta'}^n(\delta)+1)=$$
$$F_{\zeta'}(F_\zeta(\xi_\eta'')+F_{\zeta'}^n(\delta)+1)=
  F(F_{\zeta'}^n(\delta)+1)\rightarrow F(F_{\zeta'}^{n+1}(\delta))=$$
$$F_{\zeta'}(F_\zeta(\xi_\eta'')+F_{\zeta'}^{n+1}(\delta))=
  F_{\zeta'}(F_{\zeta'}^{n+1}(\delta))=F_{\zeta'}^{n+2}(\delta)$$
therefore
$$F_{\zeta'}^{n+1}(\delta)+1\rightarrow F_{\zeta'}^{n+2}(\delta)$$
Thus, in general
$$F_{\zeta'}^n(\delta)+1\rightarrow F_{\zeta'}^{n+1}(\delta)
  \hbox{ for }1\leq n<\omega\hbox{ and }F_{\zeta'}^{n+1}(\delta)<\eta_0.$$
In the case of E${}'$ ($\xi_\eta'=0$) proceed exactly as in case E,
 only replace $\delta$
by 1 and $F_\zeta(\xi_\eta'')$ by 0; then
$$F_{\zeta'}^n(1)+1\rightarrow F_{\zeta'}^{n+1}(1)
  \hbox{ for }1\leq n<\omega\hbox{ and }F_{\zeta'}^{n+1}(\delta)<\eta_0$$

Case F: One has
$\Phi_\eta(\xi+1)=F_{\xi+1}(\delta)\raeq
 \Lim_{n<\omega} F_\xi^{n+1}(F_{\xi+1}(F_\zeta(\xi_\eta''))+1)$
$$\raeq\Lim_{n<\omega} F_\xi^{1+n}(\delta)$$
so, because of property 2) b) for the sequence
 $\{\phi_{F_\xi^{1+\eta}(\delta)}\}$
$$\Phi_\eta(\xi)+1=F_\xi(\delta)+1\rightarrow\Phi_\eta(\xi+1)$$

Remark 1: In this proof the property $2\rightarrow\alpha$ (for
$2<\alpha<\eta_0$) is used.  This holds; because $\eta$ is a limit
 ordinal $<\eta_0$ and
$\eta\raeq\Lim_{x<\tau_\eta}\eta^{(x)}$,
so $\phi_\eta(1)\geq\phi_{\eta^{(1)}}(1)$, $\eta^{(0)}\geq 1$,
 $\eta^{(1)}\geq 2$;
 from which follows $2\rightarrow\alpha$.

Now that these auxiliary propositions have been proven, it is very easy
to prove 2) b) for the sequence $\{\eta_0^{(x)}\}$ of $\eta_0$:

One has an associated normal function $\Phi_{\eta_0}(\xi)$ and a limit
 ordinal $\xi_{\eta_0}$
so that
$$\eta_0=\Phi_{\eta_0}(\xi_{\eta_0})>\xi_{\eta_0}$$
and $\Phi_{\eta_0}(\xi)>\xi$ for all $\xi$ with $1\leq\xi\leq\xi_{\eta_0}$.

Since $\xi_{\eta_0}\raeq\Lim_{x<\tau_{\eta_0}}\xi_{\eta_0}^{(x)}$,
 $\eta_0^{(x)}=\Phi_{\eta_0}(\xi_{\eta_0}^{(x)})$ and by assumption
$$\xi_{\eta_0}^{(x)}+1\rightarrow\xi_{\eta_0}^{(x+1)}
  \hbox{ for }0\leq x<\tau_{\eta_0}$$
By Proposition 3, $\Phi_{\eta_0}(\xi)+1\rightarrow\Phi_{\eta_0}(\xi+1)$;
 by applying the first
Proposition (with $F(\xi)=\Phi_{\eta_0}(\xi)$) it follows that
$$\eta_0^{(x)}+1=
  \Phi_{\eta_0}(\xi_{\eta_0}^{(x)})+1\rightarrow
  \Phi_{\eta_0}(\xi_{\eta_0}^{(x)}+1)\rightarrow
  \Phi_{\eta_0}(\xi_{\eta_0}^{(x+1)})=\eta_0^{(x+1)}$$
which was to be proved.

4. Using the sequence ${\mathfrak F}_0$ of normal functions $\phi_\eta$
 of first class, whose
Existence has been proved, it is now possible to define for any limit ordinal
 $y<\phi_{F_{\Omega^{(1)}+1}}(1)$
a distinguished sequence $\{y_\eta\}$ of length $\omega$ such that
$y\daeq\Lim_{n<\omega} y_n$.

But for this we need another proposition:

Proposition 4: For all $\eta<F_\Omega(1)$ we have: if $\eta$ is
  of the second kind, $\eta\raeq\Lim_{x<\tau_\eta}\eta_x$ and
$\phi_\eta(1)=\tau_\eta$, then there exists an ordinal number $\tilde\eta$
 of the third kind with $\tilde\eta\raeq\Lim_{x<\Omega}\tilde\eta_x$
 and $\eta<\tilde\eta\leq F_\Omega(1)$,
such that $\tilde\eta_x=\eta_x$ for all $x<\tau_\eta$ and
 $\tilde\eta_{\tau_\eta}=\eta$.

Proof: The proposition holds for $\eta<\Omega^\Omega$, as we have shown
 in \S 3.  We
now suppose that the theorem holds for all $\eta'<\eta$, where $\eta$
 is a limit ordinal with
$\Omega<\eta<F_\Omega(1)$,
and show that then the theorem also holds for $\eta$:

One has $\eta=\Phi_\eta(\xi_\eta)>\xi_\eta$, for $\xi_\eta$ of the
 second kind.  You can by means of
transfinite induction readily show that
 $\phi_{\Phi_\eta(\xi)}(1a\geq\phi_\xi(1))$ for $1\leq\xi\leq\xi_\eta$,
in particular
$$\phi_\eta(1)\geq\phi_{\xi_\eta}(1)$$

If $\phi_{\xi_\eta}(1)>\tau_{\xi_\eta}$, then
 $\phi_\eta(1)>\tau_{\xi_\eta}=\tau_\eta$; so this case does not occur.

We now consider the case that $\phi_{\xi_\eta}(1)=\tau_\xi$;
 then according to the
induction hypothesis there exists a $\tilde\xi_\eta$ of third kind
 with the above required properties.

For $\eta$ then the following possibilities exist (cf.\ No.\ 2 of \S 4):

Case A: $\eta=\Omega^{\xi_0}\cdot\xi_1+\xi_2$, $\xi_\eta=\xi_2$ of the
 second kind $<\Omega^{\xi_0}$.  Then $\xi_\eta\leq\Omega^{\xi_0}$ so set
$\tilde\eta=\Omega^{\xi_0}\cdot\xi_1+\tilde\xi_\eta$.
$\tilde\xi_\eta>\Omega^{\xi_0}$, so $\xi_\eta<\Omega^{\xi_0}<\tilde\xi_\eta$
 and property 2) b) of $\tilde\xi_\eta$ is the result of
$$\phi_{\Omega^{\xi_0}}(1)>\phi_{\xi_\eta}(1);$$
since, further,
 $\phi_\eta(1)>\phi_{\Omega^{\xi_0}\cdot\xi+1}(1)
  \geq\phi_{\Omega^{\xi_0}}(1)$, $\phi_\eta(1)>\phi_{\xi_\eta}(1)$,
 therefore
$$\phi_\eta(1)>\tau_\eta$$

Case B: $\eta=\Omega^{\xi_0}\cdot(\xi_1'+1)$, $\xi_1'\geq 1$,
 $\xi_\eta=\Omega^{\xi_0}$, $\xi_0$ second kind. Then 
$$\phi_\eta(1)>\phi_{\Omega^{\xi_0}}(1)=\tau_\eta$$

Case C: $\eta=\Omega^{\xi_0}\cdot(\xi_1)$, $\xi_\eta=\xi_1$ of
 the second kind $<\Omega$.  Set $\tilde\eta=\Omega^{\xi_0+1}$.

Case D: $\eta=\Omega^{\xi_\eta''+1}$ is of third kind.

Case E and E ':
 $\eta\raeq\Lim_{\xi<\omega}F_{\zeta'}^{1+\xi}(F_\zeta(\xi_\eta'')+1)$
 or
 $\raeq\Lim_{\xi<\omega}F_{\zeta'}^{1+\xi}(1)$,
 so $\tau_\eta=\omega$,
 so certainly $\phi_\eta(1)>\tau_\eta$.

Case F: $\eta=F_{\xi_\eta}(\xi_\eta''+1)$, $\xi_\eta''\geq 1$.
Since $\eta<F_\Omega(1)$ is assumed, $\xi_\eta<\Omega$,
so $\xi_\eta=\Omega$.  Because $F_\xi(1)<\eta<F_\Omega(1)$,
$\phi_\eta(1)>\phi_{F_{\xi_\eta(1)}(1)}\geq\tau_{\xi_\eta}=\tau_\eta$.

Case G: $\eta=F_\xi(\rho+\xi_\eta)$, $0\leq\zeta<\Omega$ or $\zeta=\omega_2$;
either $\rho=F_{\zeta+1}(\mu)$, $F_{\zeta+1}(\mu)<\eta<F_{\zeta+1}(\mu+1)$,
 $\mu>0$; or then $0<\eta<F_{\zeta+1}(1)$,
and $\rho=\mu=0$.

Now if $\tilde\xi_\eta<F_{\zeta+1}(\mu+1)$, set
 $\tilde\eta=F_\zeta(\rho+\xi_\eta)$.
Then because of proposition 2
$$\tilde\eta<F_\zeta(\rho+F_{\zeta+1}(\mu+1))
  =F_\zeta(F_{\zeta+1}(\mu+1))
  =F_{\zeta+1}(\mu+1)$$
so $\tilde\eta$ satisfies the required conditions.

But if $\tilde\xi_\eta\geq F_{\zeta+1}(\mu+1)$, then
 $F_\zeta(\rho+\tilde\xi_\eta)\geq F_{\zeta_1}(\mu+1)$.
Then because $\xi_\eta<\eta<\tilde\xi$
$$\phi_\eta(1)>\phi_{\xi_\eta}(1)=\tau_\eta$$ 

In each of the above cases, $\tilde\eta$ as defined is less than
 $F_\Omega(1)$, except in
case G with $\zeta=\omega_2$; then $\eta=F_{\zeta_\eta}(1)$,
 $\tilde\eta=F_\Omega(1)$.

Remark 2: By using proposition 4 for the condition $\eta<F_\Omega(1)$
we only use an initial segment of the sequence ${\mathfrak F}_1$; because we
can for the time being forego deciding in case of F, whether
 proposition 4 for
$\eta>F_\Omega(1)$ is valid; this may be a difficult problem.

But if you only refer to Nr.\ 3 of this section
 and limit claims 1) 2) accordingly,
you can extend the sequence ${\mathfrak F}_1$ to a bigger one
with the help of further normal functions of second class,
as one can easily show;
but from a certain point we are missing again
the necessary proof.
But we do not want to do this because
we'll only use an initial segment of ${\mathfrak F}_1$.  With the help of this
initial segment, however, we get (as we shall see at once) distinguished
sequences for a significantly larger initial segment of the second-class than
Veblen.

Now, as in \S 3, we define distinguished sequences for all limit
 ordinals $y<H(1)=\phi_{F_\Omega(1)}(1)$
 (the limit $H(1)$ could of course be further advanced):

We define $\omega\daeq\ Lim_{n<\omega}(1+n)$ and assume that we have for
 each limit
ordinal $y'<y$ defined a distinguished sequence, where $\omega<y<H(1)$;
then one can also define a distinguished sequence for $y$.

One can distinguish two cases:

1) If $y$ is not in $V\phi_0$ the definition of the distinguished sequence
 is carried out very simply and exactly as in No.\ 3 of \S 3.

2) If $y$ is in $V\phi_0$,
 then there exists a largest ordinal $\eta$ with $0\leq\eta\leq F_\Omega(1)$,
 so that
$\eta$ is in $V\phi_\eta$, but not in $V\phi_{\eta'}$
for $\eta<\eta'<F_\Omega(1)$.
Again this can easily be seen with the same deduction,
as will be given in \S 5.
So we have $y=\phi_\eta(x)>x$.  $\eta$ may be in $VF_0$ or not.
Since in both cases Proposition 4 applies,
one can proceed exactly as in No.\ 3 of \S 3:

If $x$ is of the second kind, then theorem 1 of \S 2 is used.

Otherwise $x$ is of first kind.
If $\eta=0$ then
  the definition of the distinguished sequence for $y$ is obvious.
Is $\eta$ is of the first kind,
  you use theorem 2 of \S 2.
If $\eta$ is of the second kind, then theorem 3 of \S 2 is used.
This is due to $\phi_\eta(1)>\tau_\eta$;
  because if $\phi_\eta(1)=\tau_\eta$ held,
by proposition 4 a $\tilde\eta$ of third kind would exist with
 $\eta<\tilde\eta\leq F_\Omega(1)$,
  so that $y$ would be in $V\phi_{\eta+1}$;
  that would be a contradiction to the assumption
  that $y$ is not in $V\phi_{\eta'}$
  for $\eta'>\eta$.
If $\eta$ is of the third kind,
  you can also proceed exactly as in \S 3
  (Return of the case to the other cases by finitely many
   Regressions).

To conclude this section we want to compare our bound $H(1)$
to the bound $E(1)$ of the procedure of VEBLEN:

Now,
$$E(1)=\phi_{\Omega^\Omega+1}(1)=\phi_{F_0(\Omega)+1}(1)$$
$$H(1)=\phi_{F_\Omega(1)+1}=\phi_{F_{\omega_2}(\Omega)+1}(1)$$
$E(1)$ is the first critical point of $\phi_{F_0(\Omega)}$,
 $H(1)$ the first critical point of $\phi_{F_{\omega_2}(\Omega)}$.
Finally
 $$H(1)>\phi_{F_\Omega(1)}(1)=
 \phi_{F_2(1)}(1)>
 \phi_{F_1(1)}(1)>
 \phi_{F_0^2(1)}(1)=E(1).$$

\begin{center}
\S\ 5 Conditions for the possibility of a complete solution of the problem\\
 of the distinguished sequences by means of the procedure of VEBLEN
\end{center}

1. A problem of interest  would be to continue the procedure of VEBLEN
 still further and thereby to expand the sequence ${\mathfrak F}_1$
 of normal functions of first class to
a more comprehensive sequence ${\mathfrak F}$
 until you get that all the limit ordinals
 of the second 
number class have been assigned distinguished sequences.  In this
section we see what kind of properties such a sequence would have to have.

We assume that we have a sequence ${\mathfrak F}=\{\phi_\eta\}$
 of a certain length $\Lambda$
(where $\Lambda$ is a limit ordinal $\leq\omega_2$) of normal functions
 $\phi_\eta$
of first class given,
with the properties (for all $\eta<\Lambda$):
\begin{list}{}{\setlength{\leftmargin}{40pt}\setlength{\labelwidth}{35pt}
 \setlength{\itemsep}{0pt}\setlength{\topsep}{0pt}}
\item[1)]
\begin{list}{}{\setlength{\leftmargin}{15pt}\setlength{\labelwidth}{10pt}
 \setlength{\itemsep}{0pt}\setlength{\topsep}{0pt}}
\item[a)]$\phi_0(x)=\omega^x$.
\item[b)]If $\eta=\eta'+1$ then $\phi_\eta=\phi_{\eta'}'$.
\item[c)]If $\eta$ is of the second kind, $\eta$ has assigned
 a limit ordinal $\tau_\eta$ of the second number class
 and an ascending sequence $\{\eta_x\}$ of length $\tau_\eta$,
 so that $\eta\raeq\Lim_{x<\tau_\eta}\eta_x$;
 $V\phi_\eta=D_{x<\tau_\eta}V\phi_{\eta_x}$.
\item[d)]If $\eta$ is of the third kind, $\eta$ has assigned
 an ascending sequence $\{\eta_x\}$ of length $\Omega$,
 so that $\eta\raeq\Lim_{x<\Omega}\eta_x$;
 $\phi_\eta(x)=\phi_{\eta_x}(1)$.
\end{list}
\item[2)]
 For the sequences $\{\eta_x\}$ resp.\ $\{\phi_{\eta_x}\}$ in 1),
 i.e.\ for all $x$ with $1\leq x<\tau_\eta$ (in particular for
 1) d) for all $x$ with $1\leq x<\Omega$), the following hold:
\begin{list}{}{\setlength{\leftmargin}{15pt}\setlength{\labelwidth}{10pt}
 \setlength{\itemsep}{0pt}\setlength{\topsep}{0pt}}
\item[a)]
 If $x$ is a limit, then
 $\eta_x\raeq\Lim_{x'<x}\eta_{x'}$
\item[b)]$\eta_x+1\rightarrow\eta_{x+1}$
 where the arrow again has meaning as defined in \S 3, No.\ 2.
\end{list}
\item[3)]
If $\eta$ is of second kind,
$\eta\raeq\Lim_{x<\Omega}\eta_x$
and $\phi_\eta(1)=\tau_\eta$,
then there exists a certain ordinal number $\tilde\eta$
of the third kind with
$\tilde\eta=\Lim_{x<\Omega}\tilde\eta_x$
and $\eta<\tilde\eta<\Lambda$
(i.e.\ $\phi_{\tilde\eta+1}$ in ${\mathfrak F}$),
such that $\tilde\eta_x=\eta_x$ for $x<\tau_\eta$,
 $\tilde\eta_{\tau_\eta}=\eta$.
\item[4)]
If $\{\eta_\nu\}$ is a sequence (whose length is a limit ordinal $\lambda$)
 with
  $\Lim_{\nu<\lambda}=\Lambda$, then
$D_{\nu<\lambda}V\phi_\nu$ is empty.
\end{list}

In No.\ 2 of this section we will prove that you can with the help of
such a sequence ${\mathfrak F}$ solve the problem of the distinguished
 sequences really for the
whole second number class.

Conditions 1) and 2) apply to the sequences ${\mathfrak F}_0$
 (\S 3, No.\ 2) and ${\mathfrak F}_1$
(\S 4, No.\ 3).  Condition 2) needs to be proved, hence
 for the sequences $\{\phi_{\eta_x}\}$
the properties I) and II) of \S 2 are met, so that the definitions of
$\phi_\eta$ (if $\eta$ is a limit ordinal) in 1), and the theorems of \S 2
apply.  Also the further conditions may be proved, as
we will see.

2. We now want to assume the existence of a sequence
 ${\mathfrak F}$
with the above properties and draw conclusions.

Theorem 1: If $\{\eta_{x'}\}$ is any ascending sequence of length
 $\Omega$ of ordinal
numbers,
such that all $\phi_{\eta_x'}$ (for all $x<\Omega$) belongs to ${\cal F}$
(i.e., $\eta_x'<\Lambda$ for all $x<\Omega$),
then $D_{x<\Omega}V\phi_{\eta_x'}$ is empty.

Proof: Let $\eta=\Lim_{x<\Omega}\eta_x'$; then $\eta$ is of the third kind.
 If $\eta=\Lambda$, then $D_{x<\Omega}V\phi_{\eta_x'}$
is empty by property 4).

If $\eta<\Lambda$ and $\eta\raeq\Lim_{x<\Omega}\eta_x$,
then by assumption 1) d) $D_{x<\Omega}V\phi_{\eta_x}$ is empty,
since $\Lim_{x<\Omega}\phi_{\eta_x}\allowbreak(1)=\Omega$.
We are now investigating $D_{x<\Omega}V\phi_{\eta_x'}$. 

If $\alpha$ lies in $D_{x<\Omega}V\phi_{\eta_x'}$, then $\alpha$ lies in
 all $V\phi_{\eta_x'}$.  Now for any $x_0<\Omega$,
 there exists an $x_1<\Omega$, such that $\eta_{x_0}\leq\eta_{x_1'}<\eta$.

Since $\alpha$ is in $V\phi_{\eta_{x_1'}}$,
$\alpha$ is also in $V\phi_{\eta_{x_0}}$
(because $V\phi_{\eta_{x_1'}}$ is as the result of 2) for the sequence
  $\{\phi_{\eta_x}\}$ a subset of $V\phi_{\eta_{x_0}}$).

Since this holds for any $x_0<\Omega$, $\eta$ is in
 $D_{x<\Omega}V\phi_{\eta_x}$
so $D_{x<\Omega}V\phi_{\eta_x'}$
is a subset of $D_{x<\Omega}V\phi_{\eta_x}$,
so $D_{x<\Omega}V\phi_{\eta_x'}$ is the empty set.

Remark: If you drop the requirement 4) (you take
e.g.\ any portion of ${\mathfrak F}_0$ or ${\mathfrak F}_1$),
 then theorem 1 applies only if
the restriction $\Lim_{x<\Omega}\eta_x'\Lambda$ is taken as a condition.

Theorem 2: The length $\Lambda$ of ${\mathfrak F}$ is a limit ordinal of the
 third number class.

Proof: Assumption: $\Lambda=\omega_2$.
Let $\alpha$ be in $V\phi_0$.
We take in ${\mathfrak F}$ all $\phi_\eta$,
for which $\alpha$ lies in $V\phi_\eta$.
So we get a sequence $\{\phi_{\eta_\nu'}\}$ of length $\xi$
($\nu$ runs through all ordinals $<\xi$).
Either this sequence has a last element
(if $\xi$ is of the first kind),
or then $\Lim_{\nu<\xi}\eta_\nu'<\Lambda$,
because otherwise $\xi=\omega_2$,
so you would have an initial segment of the sequence $\{\eta_\nu'\}$,
whose length equals $\Omega$ and whose limit would be $<\Lambda$;
so by theorem 1 $D_{\nu<\Omega}V\phi_{\eta_\nu'}$ would be empty,
which contradicts the definition of the sequence $\{\phi_{\eta_\nu'}\}$.
Because the conclusion holds for any
  ordinal $\alpha$ from $V\phi_0$,
for every such $\alpha$ there exists an ordinal $\eta_\alpha<\omega_2$,
such that $\alpha$ is not in $V\phi_\eta$ for all $\eta$
  $\eta_0<\eta<\Lambda$.
The limit of all $\eta_\alpha$ is $H$.
Then there is no element of $V\phi_0$ in $V\phi_\eta$ for $\eta>H$.
But since $H<\Lambda$,
this yields a contradiction
that all $V\phi_\eta$ for $\eta<\Lambda$ consist of ordinals from $V\phi_0$.
-- So $\Lambda<\omega_2$.

Theorem 2': $\Lambda$ is of the third kind.

Proof: For the sake of proof we assume that $\Lambda$ is of the second kind;
then there is a sequence $\{\eta_m\}$ of length $\omega$ with
 $\Lim_{m<\omega}\eta_m=\Lambda$.
We now take a fixed $m<\omega$ and define a finite sequence
 $\{\zeta_n\}$:

Let $\zeta_0=\eta_{m+1}$.

If $\eta_m+1<\zeta_n$ and $\zeta_n=\zeta_n'+1$,
let $\zeta_{n+1}=\zeta_n'$.
Then $V\phi_{\zeta_n}$ is a subset of $V\phi_{\zeta_n+1}$.

If $\eta_m+1<\zeta_n$ and $\zeta_n$ is a limit ordinal
 with $\zeta_n\raeq\Lim_{x<\tau_{\zeta_n}}\sigma_x$,
let $\sigma_{x_0}$ be the first member of the sequence $\{\sigma_x\}$ with
  $\sigma_{x_0}>\eta_m$ and let $\zeta_{n+1}=\sigma_{x_0}$.
Then there is an ordinal $y<\Omega$,
such that the suffix of $V\phi_n$ defined by $\phi_{\zeta_n}(x)>y$
 is a subset of that
 suffix of $V\phi_{\zeta_{n+1}}$
 defined by $\phi_{\zeta_{n+1}}(x)>y$.

After finitely many steps a $\zeta_{n_0}=\eta_m+1$ is reached;
so there is an ordinal $y^{(m)}$,
such that the suffix of $V\phi_{\eta_{m+1}}$ defined by
 $\phi_{\eta_{m+1}}(x)>y$
 is a subset of that
 suffix of $V\phi_{\eta_{m+1}}$
 defined by $\phi_{\eta_{m+1}}(x)>y$;
the considered suffix of $V\phi_{\eta_{m+1}}$ consists only of critical
 points of $V\phi_{\eta_m}$.

For every $m<\omega$ there exists such a $y^{(m)}$;
denote the limit of all $y^{(m)}$ as $y_0$.
Then $y_0<\Omega$,
and the intersection of the suffixes of $V\phi_{\eta_m}$ defined by
 $\phi_{\eta_m}(x)>y_0$ is not empty
(cf.\ requirement I) of \S 2),
so $D_{m<\omega}V\phi_{\eta_m}$ is not empty,
which contradicts requirement 4) of this section.

Theorem 3:
If you have uniquely constructed a sequence of length $\Lambda$, of the
 third kind $<\omega_2$
then the problem of the distinguished sequences for the whole
second number class is solvable.

Proof: a) We set $\omega\daeq\Lim_{n<\omega}(1+n)$ and
  suppose,
we have defined a distinguished sequences for all the limit ordinals $y'<y$,
where $y$ is a limit ordinal with $\omega<y<\Omega$ and
show that a distinguished sequence for $y$
  can be defined.

b) If $y$ is not in $V\phi_0$, so
  $y=\omega^{x_0}\cdot x_1+x_2$,
where $1\leq x_0<\omega$, $1\leq x_1<\omega$,
  $0\leq x_2<\omega^{x_0}$.
you can see that the distinguished sequence of $y$ is just like that
  defined in Nr,\ 3 of \S 3.

c) If $y$ lies in $V\phi_0$,
then there is a last normal function $\phi_\eta$ of ${\mathfrak F}$,
such that $y$ is in $V\phi_\eta$,
but not in $V\phi_{\eta'}$ for all $\eta'$ with $\eta<\eta'<\Lambda$.
This can be shown as follows:

Let $\{\eta_\nu'\}$ be the sequence of length $\lambda$,
such that $y$ is in all $V\phi_{\eta\nu'}$,
and only in these sets of values
(where $\nu$ goes through all the ordinals $<\lambda$).
We assume
this sequence has no last element,
i.e.\ $\lambda$ is a limit ordinal.
Set $\eta'=\Lim_{\nu<\lambda}\eta_\nu'$,
from 4) immediately follows $\eta'<\Lambda$.
Were $\eta'$ of the third kind,
you could show exactly as in the proof of theorem 1,
that $D_{\nu<\lambda}V\phi_{\eta_\nu'}$ would be empty;
So you would have a contradiction.
Were $\eta'$ of the second kind,
so $\eta'\raeq\Lim_{x<tau_{\eta'}}\eta_x''$,
one could show that, as in the proof of theorem 1,
$D_{\nu<\lambda}V\phi_{\eta'}$ is a subset of
$D_{x<\tau_{\eta'}}V\phi_{\eta_x''}=V\phi_{\eta'}$,
i.e. that $y$ would also be in $V\phi_{\eta'}$,
which contradicts the requirement.
So the sequence $\{\eta_\nu'\}$ has a final element,
that we call $\eta$.
- So you have $y=\phi_\eta(x)>x$.

d) One can now define a distinguished sequence for $y$ as in the
  previous section
(cf. \S 3, No. 3).
We summarize only briefly:

Let $x$ be first kind.
Then $\eta=0$,
so the definition of the sequence of $y$ is obvious.
Is $\eta$ first kind,
so use theorem 2 of \S 2.
Is $\eta$ of the second kind,
$\eta\raeq\Lim_{x<\tau_\eta}\eta_x$,
then use theorem 3 of \S 2 for the sequence $\{\phi_{\eta_x}\}$.
This last definition fails,
if $y=\phi_\eta(1)=\tau_\eta$.
But that can not happen
otherwise, assuming condition 3), an ordinal $\tilde\eta>\eta$ would exist,
so $y$ would be in $V\phi_{\eta+1}$.
If $\eta$ is of third kind,
then the case can be reduced through finitely many regressions to the other
 cases.
As shown in \S 3,
in the case of the application of theorem 3 of \S 2 no failure of the
  definition can occur.

If $x$ is of the second kind, then use the first theorem of \S 2.

At the proofs of the restricted theorem 1 and the theorem 2 we have
 that only conditions 1) and 2) are used;
Conditions 3) and 4) are only essential for Theorem 2' and 3.
Theorem 2 says
that the sequence ${\mathfrak F}$ while preserving their properties 1) and 2)
 can not continue through the entire third class of ordinals,
as it may seem at first.
The third theorem is the problem of the distinguished sequences on the
(not less difficult)
of existence
(or the construction)
from the preceding.
We can clearly define
what we mean by continuation of the procedure from VEBLEN to the
 desired complete solution of the problem of distinguished
 sequences:
below it is supposed to form a sequence of normal functions with the above
 given conditions.
There is no sequence of normal functions with these properties
(which is possible),
so you can not solve the problem of distinguished sequences with the method
 of VEBLEN.

Remark: The train of thought of \S\S 3 to 5 can by the way
  be easily generalize to higher number classes or normal functions of
  higher classes.
The problem of the distinguished sequences would then be generalized:
For every limit ordinal $\eta<\omega_{k+1}$
(where $k$ is any fixed integer)
choose a unique ascending sequence $\{\eta_x\}$,
whose length $\lambda$ is a limit ordinal with $\omega\leq\lambda\leq\omega_k$,
such that $\eta=\Lim_{x<\lambda}\eta_x$.

\begin{center}
\S\ 6 The analog for number-theoretic functions
\end{center}

1.  For the monotonically increasing number-theoretic functions $f(x)$
(as a normal functions of 0th class)
can be given analogous considerations.
Here the argument goes through all natural numbers and
 $\Lim_{x<\omega}f(1+x)=\omega$.
Since such a function does not have critical points,
let us substitute another operation for the derivation,
which results in a dilution of the values of the function,
e.g.\ the formation of the iteration $g(x)=f^2(x)=f(f(x))$ from $f(x)$.
If one uses only such number-theoretic functions $f(x)$,
which satisfy the condition $f(1)>1$
(which we also assumed in the normal functions of first class),
then $g(1)>f(1)$ and $Vg$ is a subset of $Vf$.
Similar to the method of VEBLEN one can then for a given
 initial function,
e.g. $2^x$,
form their iteration,
from the new function again the iteration, etc.,
then take the first values of the generated function sequence of type $\omega$
 etc.

The sequence of normal functions of first class in \S 5 would thus become a
  sequence ${\mathfrak F}'=\{f_n\}$ of length $\lambda$
(where $\lambda$ is a limit ordinal $\leq\Omega$)
of number-theoretic functions $f_\eta$,
with the properties
(valid for all $\eta<\lambda$):

1) $f_0(x)=2^x$ (e.g.).

2) If $\eta=\eta'+1$, then $f_\eta=f_{\eta'}^2$.

3) If $\eta$ is of the second kind,
so $\eta$ has an associated ascending sequence $\{\eta_x\}$ of length $\omega$
so that $\eta\raeq\Lim_{x<\omega}\eta_x$,
then $\eta_x+1\rightarrow\eta_{x+1}$ for all $x<\omega$,
and $f_\eta(x)=f_{\eta_x}(1)$.
We introduce the abbreviation $\alpha\rightarrow\beta$ for:
$Vf_\eta$ is a subset of $Vf_\alpha$ for all $\eta$ with
  $\alpha<\eta\leq\beta$.

The other conditions
which we set up for ${\mathfrak F}$,
are invalid or not of interest here.

2. Assuming the existence of such a sequence ${\mathfrak F}'$ of
  number-theoretic functions one deduces according to \S 5 the theorems:

Theorem 1: For each sequence $\{f_{\eta'}\}$ of length $\omega$ of
  number-theoretic functions from ${\mathfrak F}'$
  $D_{n<\omega}Vf_{\eta_n'}$ is empty,
provided $\Lim_{n<\omega}\eta_n'\neq\lambda$\footnote{Whether the theorem
also applies in the case $\Lim_{n<\omega}\eta_n'=\lambda$ in general,
is questionable.}.

Proof: Let $\eta=\Lim_{n<\omega}\eta_n'$.
Because $\eta<\lambda$ $f_\eta$ exists,
and one has a sequence $\eta\raeq\Lim_{x<\omega}\eta_x$.
Exactly as in the proof of theorem 1 of \S 5,
$D_{n<\omega}Vf_{\eta_n'}$ is a subset of the empty set
  $D_{x<\omega}Vf_{\eta_x}$,
so is empty.

Theorem 2: $\lambda<\Omega$.

Proof: Assumption: $\lambda=\Omega$.
Let $\alpha$ be an element of $Vf_0$.
We take all $f_\eta$ in ${\mathfrak F}'$,
with $\alpha$ in $Vf_\eta$;
so we get an ascending sequence $\{\eta_\nu'\}$ and
 a sequence $\{f_{\eta_\nu}'\}$ of number-theoretic functions, of
  length $\xi$.
Either this sequence has a final element
(if $\xi$ is of the first kind),
or $\Lim_{\nu<\xi}\eta_\nu'<\lambda$
(because otherwise $\xi=\Omega$ would hold,
so a subsequence of length $\omega$ of the sequence $\{\eta_\nu'\}$
 would exist,
whose limit is less than $\lambda$,
which would immediately lead to a contradiction because of theorem 1).

Thus, for every $\alpha$ from $Vf_0$ there exists an ordinal $\eta_0<\Omega$,
so $\alpha$ is not in $Vf_\eta$ for all $\eta$ with
$\eta_\alpha<\eta<\lambda$.  Let $H$ be the limit of the $\eta_\alpha$,
so $H<\Omega$ and no $\alpha$ from $Vf_0$ is in $Vf_\eta$ for
 $\eta>H$,
which is impossible.
So the assumption $\lambda=\Omega$ is to be rejected.

A sequence with the properties of ${\mathfrak F}'$ can not be of
 length $\Omega$
(by such an assignment of number-theoretic functions to the
ordinal numbers an uncountable a priori
well-order in the continuum would thereby be attainable; cf. No. 2 of § 1)

3. Application.
It is not possible to assign to
all of the limit ordinals $\eta$ of the third class sequences $\{\eta_x\}$
  of lengths $\tau_\eta<\Omega$,
so that the conditions 2) a) and 2) b) of No.\ 2 of \S 3 are met.
Otherwise you could conclude the existence of a sequence of length $\omega_2$
  of normal functions of first class having the conditions 1) 2) of No.\ 1
  from \S 5,
which is impossible according to \S 5
(because exactly like in No.\ 2 of \S 3 one would exclude  from 2) b)
of \S 3 the condition 3) of \S 3,
i.e.\ one can prove the condition 2) b) of \S 5 and thus also 1) of \S 5
  for all $\eta<\omega_2$).

The conditions 2) a) and 2) b) of \S 3 correspond in the
  second number class to the following problem of distinguished sequences
  with a constraint:
Each limit ordinal $\eta$ of the second number class has a unique ascending
  sequence $\{\eta_n\}$ of length $\omega$ assigned,
so that $\eta\raeq\Lim_{n<\omega}\eta_n$,
and for $0\leq\eta<\omega$,
$\eta_n<\beta\leq\eta_{n+1}$, and
$\beta\raeq\Lim_{n<\omega}\beta_n$
 the relation $\beta_0\geq\eta_n$ follows.

It is very easy to see
that solutions to this problem exist,
if you limit yourself to any initial segment of the second number class
(however, you can only select a clear solution if you have clearly solved
 the problem of the distinguished sequences without constraint).

Set $\omega\raeq\Lim_{n<\omega}(1+n)$.
We assume
all limit ordinals $\eta$,
which is smaller than a certain limit ordinal $x>\omega$ of the second
 number class,
have associated sequences
$\{\eta_n\}$ with the above constraint and show
that all limit ordinals $\leq x$ can then be assign such sequences:

Let $\{x_\eta\}$ be an ascending sequence for $x$;
if $\eta$ is a limit with $x_n<\eta\leq x_{n+1}$
($0\leq n<\omega$)
and associated sequence $\{\eta_n\}$,
replace the sequence $\{\eta_n\}$ by the suffix,
whose first member is $\geq x_n$;
if $\eta\leq x_0$,
let the sequence $\{\eta_n\}$ be changed.
The new sequences with the sequence $\{\eta_n\}$ satisfy the constraint.

Since with this induction proof with each step the
 already existing sequences is changed,
you can not conclude
that the problem is solvable for the whole second number class.
It even holds the

Theorem: With the above constraint, the problem of the distinguished
 sequences
can not be solved for all limit ordinals of the second number class.

Proof: We show that assuming the opposite, that a sequence
of number-theoretic functions with the properties 1) to 3) of
 ${\mathfrak F}'$
(where the sequence $\{\eta_x\}$ in condition 3) is just the
 given sequence of $\eta$ with the above constraint),
which would be of length $\Omega$,
could be constructed according to the remarks of this section,
 is impossible:

Set $f_0(x)=2^x$.
If all $f_{\eta'}$ for $\eta'<\eta$
(where $1\leq\eta<\Omega$)
have been defined  satisfying the required conditions,
you can also define $f_\eta$:

If $\eta=\eta'+1$, then $f_\eta=f_{\eta'}^2$.

If $\eta$ is of the second kind and $\{\eta_x\}$ the given sequence
  for $y$,
so for fixed $m<\omega$ we define a finite sequence $\{\zeta_n\}$
  through the stipulations:
Let $\zeta_0=\eta_{m+1}$.
If $\eta_m+1<\zeta_n$ and $\zeta_n=\zeta_n+1$,
set $\zeta_{n+1}=\zeta_n'$;
so $\eta_{m+1}\leq\zeta_{n+1}<\zeta_n$ and
  $\zeta_{n+1}\rightarrow\zeta_n$.
If $\eta_m+1<\zeta_n$,
$\zeta_n$ is of second kind, and $\{\zeta_n^{(x)}$ is the
  ascending sequence for $\zeta_n$,
let $\zeta_{n+1}=\zeta_n^{(1)}$;
then $\eta_m+1\leq\zeta_{n+1}<\zeta_n$ and
  $\zeta_{n+1}\rightarrow\zeta_n$.
After finitely many regressions, $\zeta_{n_0}=\eta_m+1$;
this is followed by $\eta_m+1\rightarrow\eta_{m+1}$,
and you can set $f_\eta(x)=f_{\eta_x}(1)$.


\begin{center}
\S\ 7 Reduction of the problem of the distinguished sequences\\
 to other problems related to normal functions
\end{center}

1. As you can easily show
(see note 2 of this paragraph),
you can have distinguished sequences for all limit ordinals of the second
 number class determined,
if you have such already available for all ordinals of $V\phi_0$,
i.e.\ for all ordinals $\omega^x$ ($1\leq x<\Omega$).
Therefore, we want to limit ourselves to $V\phi_0$ in this appendix.
Furthermore, we introduce the following abbreviation: If
$\{y_\eta\}$ is an ascending sequence of length $\omega$ of ordinal numbers,
let $\Delta y_n=y_n-y_{n-1}$
(defines for $1\leq n<\omega$) be the corresponding difference sequence.

We now prove the equivalence of the following four problems:

\begin{list}{}{\setlength{\leftmargin}{40pt}\setlength{\labelwidth}{35pt}
 \setlength{\itemsep}{0pt}\setlength{\topsep}{0pt}}
\item[1)]The problem of the distinguished sequences for the set $V\phi_0 $;
i.e.\ any ordinal $y$ of $V\phi_0$ has assigned a distinguished sequence
 $0<y_0<y_1<y_2<\dots$ of length $\omega$ with
 $\Lim_{n<\omega}y_n=y$.
\item[2)]The above problem with the following constraints:
\begin{list}{}{\setlength{\leftmargin}{15pt}\setlength{\labelwidth}{10pt}
 \setlength{\itemsep}{0pt}\setlength{\topsep}{0pt}}
\item[a)]$y_0<\Delta y_1$;
\item[b)]$\Delta y_n\leq\Delta y_{n+1}$ for $1\leq n\leq\omega$.
\end{list}
\item[3)]Any number $y$ of $V\phi_0$ has assigned a normal function $\psi_y$
 of first class,
such that $y$ is the first critical point of $\psi_y$ ($y=\psi_{y'}(1)$.
\item[4)]
To any given normal function $f(x)$ of first class,
whose values are only ordinals of $V\phi_0$,
to uniquely assign a normal function $F(x)$ such that $F'(x)=f(x)$.
\end{list}

For the equivalence of the problems the following property 
 must be proved:
If you have a solution to such a problem, you can also for the
 other problems above construct unique solutions.

Proof of equivalence:

1) We first prove the equivalence of the first two problems.
We just have to show that we can from a clear solution of the first
problem construct a unique solution to the second problem
(because the converse is trivial):

Suppose $y=\omega^x$ and $x=x'+1$.
Then set $y_n=\omega^{x'}\cdot(1+n)$,
thus one has a sequence $\{y_n\}$,
that meets the constraints of the second problem.

But if $y=\omega^x$,
$x$ is a limit and $\{y_n\}$ is the
  given distinguished sequence for $y$ satisfying the condition,
then for every $n<\omega$ we can choose the smallest ordinal $x_n>0$
such that $y_n\leq\omega^{x_n}$.
Then $\Lim_{n<\omega}\omega^{x_n}=y$.
In this sequence, members appear several times,
so you take them only once.
The resulting ascending subsequence $\{z_n\}$ satisfies
  conditions a) and b);
because all the $z_n$ are $V\phi_0$,
and after the second proposition of \S 4
$$\Delta z_n=z_n\hbox{ for }1\leq n<\omega$$

2) Next we prove the equivalence of the second and third problem:

a) Suppose one has a given solution to the second problem.
Then you can solve the third problem clearly:
Let $y$ be an arbitrary ordinal from $V\phi_0$.
By assumption, $y$ has a sequence $\{y_n\}$ with the
  constraints a) and b) satisfied.
We now define the normal function $\psi_y$ as follows:

Let $\psi_y(x)=y_0+x$ for $1\leq x<y_0$,
further
$$\psi_y(y_0)=y_0\cdot2$$

For given $n$
$$\psi_n(y_n+x)=y_{n+1}+x\hbox{ for }1\leq x<\Delta y_0$$
$$\psi_y(y_{n+1})=y_{n+1}+\Delta y_{n+1};$$
this then holds for $0\leq n<\omega$.

Finally, let $\psi_y(x)=x$ for all $x$ with $y\leq x<\Omega$.

$\psi_y(x)$ is a normal function;
because first, it is monotonically rising
(because of $\Delta y_n\leq\Delta y_{n+1}$)
and second, it is continuous:
if $y_{n+1}$ is a limit ordinal,
then
 $\Lim_{x<\Delta y_{n+1}}\psi_y(y_n+x)=y_{n+1}+\Delta y_{n+1}=
  \psi_y(y_{n+1})=\psi_y(\Lim_{x<\Delta y_{n+1}}(y_n+x))$
  for $0\leq n<\omega$,
and if $y_0$ is a limit ordinal,
then
 $\Lim_{x<y_0}\psi_y(x)=y_0+y_0=\psi_y(y_0)=
  \psi_y(\Lim_{x<y_0}(1+x))$;
Furthermore, because of
 $y_n<\psi_y(y_n)\leq y_n+\Delta y_{n+1}=y_{n+1}$,
 $\Lim_{n<\omega}\psi_y(y_n)=\Lim_{x<y}\psi_y(1+x)=y$.
In addition, $y$ is the first critical point;
because $1\leq y'<y$, you have three cases:
$y'<y_0$, $y'=y_n$ (for a given $n$),
or $y'=y_n+x$ (for a given $n$ and $1\leq x<\Delta y_{n+1}$1).
In the first case, $\psi_y(y')>y'$ because $y_n>0$;
in the second case $\psi_y(y')=y_n+\Delta y_n>y_n=y'$,
because $\Delta y_n>0$;
in the third case would be $\psi_y(y')=y'$,
that is, $y_{n+1}+x=y_n+x$,
so $\Delta y_{n+1}+x=x$;
but because the first critical point of the normal function
 $\chi(x)=\Delta y_{n+1}+x$
 is $\Delta y_{n+1}\cdot\omega$
$x\geq\Delta y_{n+1}\cdot\omega$ would hold,
which is a contradiction to the assumption $x<\Delta y_{n+1}$.

b) Converse: Suppose the third problem is clearly solved
and $y$ is an element of $V\phi_0$;
by assumption, one has the normal function $\phi_y(x)$.
Then, by Theorem 2 of \S 2, $y=\Lim_{n<\omega}\psi_y^n(1)$.
$\{\phi_y^n(1)\}$ is an ascending sequence
because 1 is not a critical point of $\psi_y$.
Further, (see No. 3 of \S 1)
$$\psi_y^0(1)=1\leq\psi_y(1)\hbox{---}1\leq\psi_y^2(1)\hbox{---}\psi_y(1)\leq$$
$$\psi_y^3(1)\hbox{---}\psi_y^2(1)\leq\dots$$
i.e.\ our sequence for $y$ satisfies the constraints a) and b).

3) Finally, we prove the equivalence of the third and fourth problems:

a) A clear solution of the third problem is followed by a clear solution
 of the fourth:
Let $f(x)$ be a given normal function whose
 value set is a subset of $V\phi_0$.
Now set
$$F(f(x)=f(x)$$
$$F(f(x+y)=f(x)+\psi_{f(x+1)}(y)\hbox{ for }1\leq y<f(x+1);$$
this then holds for $1\leq x<\Omega$.

Our function $F$ is a normal function,
as can be demonstrated by using Proposition 2 of \S 4;
furthermore $F(f(x))=f(x)$,
$F(y)>y$ for $y\neq f(x)$.

b) Converse: The fourth problem has clearly been solved,
and $y$ is an arbitrary ordinal from $V\phi_0$;
then we define the normal function $f(x)$ by specifying:
$Vf$ consists of exactly all ordinals of $V\phi_0$,
which are $\geq y$.
By assumption, there exists a well-defined normal function $F(x)$,
such that $F'(x)=f(x)$, so $y=F'(1)$.
Set $\psi_y(x)=F(x)$,
so you have a solution to the third problem.

2. Comments:

1) The second problem is a problem of distinguished
 sequences with a constraint.
If you limit yourself to $V\phi_0$,
it has solutions.
But it is not solvable for all limit ordinals of the second number class;
because e.g.\ $y=\epsilon+\omega$ and $y=\Lim_{n<\omega}y_n$,
where $\{y_n\}$ is any ascending sequence for $y$.
Then let $y_{\eta_0}$ be the first member of the sequence with
 $y_{n_0}\geq\epsilon$.
If $n_0=0$, then you have $y_0>\Delta y_1$;
if $n_0>0$, then $\Delta y_{n_0}\geq\epsilon$,
 $\{Delta y_{n_0+1}<\epsilon$,
that is $\Delta y_{n_0}>\Delta y_{n_0+1}$.
In both cases, then, the constraints of the second problem are
 not fulfilled.

On the other hand, if the constraint a) is dropped,
you have the problem of distinguished sequences with only one
 insignificant secondary condition b);
this is solvable for the whole second number class;
because you can easily find the equivalence of this problem with the problem
 of proving the existence of distinguished sequences without constraint:

You just have to show
that from a clear solution of the latter a clear solution
 of the former can be constructed
(the converse is trivial):
For $\omega$, $\omega=\Lim_{n<\omega}(1+n)$.
Let $y>\omega$ be a limit ordinal of the second number class with
 a given distinguished sequence,
and we have constructed sequences for all the limit ordinals $y'<y$,
satisfying b).
Then one can also construct such a sequence for $y$:
If $y$ is in $V\phi_0$,
you can even construct a sequence
complying with a) and b)
(according to the above argument in No. 1 of this section).
If $y$ is not in $V\phi_0$,
we write $y=\omega^{x_0}\cdot x_1+x_2$,
with $1\leq x_0<\omega$, $1\leq x_1<\omega$, $0\leq x_2<\omega^{x_0}$.

Then $x_2$ is of  second kind and $\{\zeta_n\}$ is the corresponding sequence
  with the property b),
$\omega^{x_0}\cdot x_1+\zeta_n$ is such a sequence for $y$.
If $x_2=0$, $x_1=x_1'+1$($x_1'\geq 1$)
and $\{\zeta_n\}$ is the sequence with property b)
belonging to $\omega^{x_0}$,
then $\{\omega^{x_0}\cdot x_1'+\zeta_n\}$ is
such a sequence for $y$.

2) The essential result of this section, however, is the
  reduction of the problem of the distinguished sequences
  to the third or fourth problem;
these problems (for $V\phi_0$) are equivalent to the problem of
  distinguished sequences for all limit ordinals of the second class,
as one can easily show;
because of the proven equivalence of the four problems, one only has to show
that the first problem is equivalent to the problem of the
 distinguished sequences
(for all limit ordinals of the second number class):

We assume that the problem of the distinguished sequences is solved for all
  ordinals from $V\phi_0$ and for all limit ordinals $<y$,
where $y$ is a limit with $\omega<y<\Omega$,
that is not in $V\phi_0$.
Then you can define an distinguished sequence for $y$,
by writing $y=\omega^{x_0}\cdot x_1+x_2$
  and takes the same definitions as in No. 3 of \S 3.

\begin{center} References \end{center}

(1) G.\ H.\ HARDY, Quarterly Journal of Mathematics, vol.\ 35 (1903).

(2) 0.\ VEBLEN, Transactions of the American Mathematical Society, vol.\ 9
 (1908).

(3) A.\ CHURCH und S.\ C.\ KLEENE, Fundamenta Mathematicae, Bd.\ 28 (1937).

(4) F.\ HAUSDORFF, Mengenlehre, Dover Publications, 3rd Revised Edition
 (New York 1944).

(5) G.\ CANTOR, Gesammelte Abhandlungen (Berlin 1932).

(6) G.\ CANTOR, Mathematische Annalen, Bd.\ 49 (1897).

\end{document}